\documentclass{article}

\usepackage {amssymb, amsmath}

\begin{document}

\def\SL{{\rm SL}}
\def\GL{{\rm GL}}
\def\PGL{{\rm PGL}}
\def\PO{{\rm PO}}
\def\O{{\rm O}}
\def\SO{{\rm SO}}

\def\R{{\Bbb R}}
\def\C{{\Bbb C}}
\def\S{{\Bbb S}}
\def\E{{\Bbb E}}
\def\PE{{\Bbb P\Bbb E}}
\def\OO{{\Bbb O}}
\def\P{{\Bbb P}}
\def\K{{\frak K}}
\def\W{{\cal W}}
\def\F{{\cal F}}

\def\bJ{{\bold J}}
\def\bI{{\bold I}}
\def\bK{{\bold K}}
\def\bM{{\bold M}}

\def\fra{{\frak a}}
\def\frb{{\frak b}}
\def\frc{{\frak c}}
\def\frd{{\frak d}}
\def\fre{{\frak e}}
\def\frf{{\frak f}}
\def\frg{{\frak g}}
\def\frh{{\frak h}}

\def\frl{{\frak l}}
\def\frm{{\frak m}}
\def\frn{{\frak n}}
\def\frp{{\frak p}}
\def\frq{{\frak q}}
\def\frr{{\frak r}}
\def\frs{{\frak s}}
\def\frt{{\frak t}}
\def\fru{{\frak u}}
\def\frv{{\frak v}}
\def\frw{{\frak w}}
\def\frx{{\frak x}}
\def\fry{{\frak y}}
\def\frz{{\frak z}}

\def\le{\leqslant}
\def\ge{\geqslant}

\def\frA{{\frak A}}
\def\frB{{\frak B}}
\def\frC{{\frak C}}

\def\Op{{\rm Op}}

\def\tr{{\rm tr \,}}
\def\codim{{\rm codim\,}}

\def\emptyset{\varnothing}
\def\phi{\varphi}
\def\epsilon{\varepsilon}

\def\nul{{nul}}
\def\fin{{fin}}
\def\dist{{\rm dist}}
\def\sol{{\rm solv}}
\def\far{{\rm far}}
\def\Pen{{\rm Pen}}
\def\Ass{{\rm Ass}}
\def\Karp{{\rm Karp}}
\def\Pass{{\rm Pass}}
\def\Flat{{\rm Str}}
\def\karp{{\rm karp}}
\def\pass{{\rm pass}}

\def\gen{{\rm gen}}

\newcounter{sec}
 \renewcommand{\theequation}{\arabic{sec}.\arabic{equation}}

\begin{center}
\Large\bf

Pencils of geodesics in symmetric spaces,

Karpelevich boundary,

and associahedron-like polyhedra

\smallskip

\large\sc

Yurii A. Neretin

\end{center}

\begin{flushright}

\it
To the memory of F.I.Karpelevich

\end{flushright}

These notes on  the matrix geometry are prepared for
the Karpelevich
memorial volume of AMS. Our standpoint is
the geometric part
  of his  treatise
'The geometry of geodesics and the eigenfunctions
of the Beltrami--Laplace operator
 on  symmetric spaces' (1966).
The subject of analytical part of his work
was the Dynkin's problem of description
of  the Martin boundary for symmetric
spaces, for its solution see \cite{Olsh}, \cite{GJT2}.
We do not touch this subject.

The existence of
the complicated Karpelevich
boundary is well known, but   in few
works (I know only \cite{Kush}, \cite{GJT1}, \cite{GJT2})
 it was really discussed.
We give elementary geometric descriptions
of the Karpelevich boundary and
of some Karpelevich-like constructions
We consider only the spaces
$\GL(n,\R)/\O(n)$ and use a minimal necessary language.

  Boundaries of symmetric spaces are
an old  subject arising
to the works of Chasles (1864--65),
Schubert (1879),  Study (1886),
and Semple (1946--52) on the enumerative algebraic
geometry. Later these boundaries
 appeared as objects and tools
of the analysis
on symmetric spaces. Some references are
\cite{Sem1}-\cite{Sem3}, \cite{Sat}, \cite{OS},
\cite{DCP}, \cite{Pop}, \cite{Olsh},
 \cite{Ner2}, \cite{Ner4}, \cite{Ner7},
\cite{GJT2}.
For further references and for the history of the subject,
see \cite{Ner4}, \cite{GJT2}, \cite{Kle}.

In the present work we start from
 pencils of geodesics
and gluing of points at infinity as limits
of pencils, i.e., we begin from the ordinary
differential
geometry. As a result, we obtain some elements
of the geometry of angles at infinity.

Recall, that in symmetric spaces  the usual distance
 is replaced by the so-called
``complex distance'' (other terms are
``compound distance'', ``composite distance'',
``angles'', and ``stationary angles'').
 This ``distance'' is a
finite collection of real numbers.
Not much is known about this geometrical structure.
In last years, after \cite{Klya}, the problem
 of the triangle inequality
became popular (see \cite{Ful}, \cite{Ner6}, \cite{KT}, \cite{KLM}).
Another fact of the geometry of angles is the
``compression
of angles'' phenomenon (see \cite{Ner1},
 \cite{Ner3}, VI.3, \cite{Ner5}, \cite{Kou}).

We show, that  the 'geometry
of angles at infinity' leads to some moduli
space
like polyhedra as the associahedron,
the permutoassociahedron,
and the Karpelevich polyhedron;
 the associahedron was
constructed by Stasheff  \cite{Sta},
see also  \cite{DJS}, \cite{Kapo},
 it is a real form of the Deligne--Mumford
moduli space of point configurations on
a rational curve, the permutoassociahedron
was constructed by Kapranov  \cite{Kap}
(see also \cite{RZ}, \cite{DCP2}),
an algebraic-geometric counterpart
of the Karpelevich polyhedron (polydiagonal blowing)
 was recently constructed
 by Ulyanov \cite{Ulya}).

Our Section 1 contains preliminaries on the symmetric spaces
$\GL(n,\R)/\O(n)$.

In Section 2, we give an explicit description
of the Satake--Furstenberg boundary of these
spaces.

In Section 3, we discuss pencils of geodesics
in symmetric spaces. Following Karpelevich, we define
finite pencils, null pencils and solvable pencils.

After this, we add limits of pencils as points of
a symmetric space
at infinity. Three types of pencils give  4 different boundaries.

The most simple case is discussed in Section 4,
 limits of all the finite pencils form the so-called visibility
sphere
at infinity.

Limits of solvable pencils form a noncompact
space. This space can be compactified by two
similar inductive procedures (Section 5).
In these cases, the compactifications of  Cartan (flat)
subspaces  are some combinatorial
 polyhedrons, namely the permutoassociahedrons
and the Karpelevich polyhedrons.
They are described explicitly in Section 6.
In Section 7, we finish the description
of the associahedral and Karpelevich
boundaries.

In Section 8, we construct
the sea urchin \cite{Ner7} using  null pencils.

\medskip

{\bf 1.
Symmetric spaces $\GL(n,\R)/\O(n)$ and $\PGL(n,\R)/\PO(n)$
}

  \addtocounter{sec}{1}
 \setcounter{equation}{0}

\medskip

{\bf 1.1. The space $\GL(n,\R)/\O(n)$.}
Consider the space $\E_n$ of positive
 definite real symmetric matrices of a size
$n\times n$.
The general linear group $\GL(n,\R)$
acts on this space by the transformations
$
X\mapsto gXg^\top
$,
where $X\in\E_n$, $g\in\GL(n,\R)$;
and the sign $\top$ denotes the transposition.
The stabilizer of the point $X=E$ is the orthogonal
 group
$\O(n)$ and hence
$
\E_n=\GL(n,\R)/\O(n).$

{\bf 1.2. The space $\PGL(n,\R)/\PO(n)$.}
Denote by $\PE_n$ the space of positive definite
$n\times n$ matrices defined up to a scalar factor:
$$X\sim \lambda X, \qquad \text{where\quad
$X\in \E_n$, $\lambda>0$}.$$

Obviously,
$\PE_n=\PGL(n,\R)/\PO(n)$;
where $\PGL(n,\R)$
is the quotient group of $\GL(n,\R)$
 by the  group $\R^*$ of scalar matrices, and
 $\PO(n)=\O(n)/\{\pm 1\}$.

We also can consider the space $\PE_n$ as the space of
 all positive definite matrices $X$ such that $\det(X)=1$,
hence $\PE_n=\SL(n,\R)/\SO(n).$

{\bf 1.3. Quadratic forms.}
For each $X\in \E_n$ we define the positive definite
bilinear form on $\R^n$ by
\begin{equation}
Q_X(v,w)=
\tfrac 12 \sum\nolimits_{i,j\le n} x_{ij}v_i w_j
,
\end{equation}
where $v=(v_1,\dots,v_n)$, $w=(w_1,\dots,w_n)\in\R^n$,
and $x_{ij}$ are the matrix elements of the matrix $X$.

Thus we identify $\E_n$ with the space of positive definite quadratic forms
and $\PE_n$ with the space of positive
 definite quadratic forms defined
up to a scalar factor.

{\bf 1.4. Space of ellipsoids.}
For any form (1.1) we consider the ellipsoid
\begin{equation}
\tfrac 12 \sum\nolimits_{i,j\le n} x_{ij}v_i v_j=1
.\end{equation}
Thus we identify the space $\E_n$
with the space of ellipsoids with center
at 0.
Also we identify the space $\PE_n$
with the space of ellipsoids
defined up to a homothety
$
v\mapsto \lambda v$, where
 $v\in\R^n$, $\lambda>0$.

{\bf 1.5. Complex distance.}
Let $X$, $Y\in \E_n$.
We consider the equation
$$
\det(X-\lambda Y)=0
$$
and denote its solutions
(they are real and positive)   by
$$
\lambda_1(X,Y)\ge
\lambda_2(X,Y)\ge
\dots\ge
\lambda_n(X,Y)
.$$

{\sc Theorem  1.1.}
{\it For $X$, $Y$, $X'$, $Y'$
the following conditions are equivalent.

(i) There exists $g\in\GL(n,\R)$ such that
$ gXg^\top=X'; \qquad gYg^\top =Y'$.

 (ii) $\lambda_j(X,Y)=\lambda_j(X',Y')$
for all $j$.}

We define the {\it complex distance} in $\E_n$
as the collection
$$
\psi_j(X,Y)=\ln \lambda_j(X,Y)
.$$

In the space $\PE_n$, this collection is defined
up to a common additive constant
$$
(\psi_1,\dots,\psi_n)\sim
(\psi_1+\tau,\dots,\psi_n+\tau)
.$$


{\sc Theorem 1.2.} (\cite{Ner6})
{\it Fix $X$, $Y$, $Z\in\E_n$.  Let
$\Psi=(\psi_1,\dots,\psi_n)$ be the complex distance between
$X$ and $Y$,
$\Phi=(\phi_1,\dots,\phi_n)$ be the complex distance between
$Y$ and $Z$,
$\Theta=(\theta_1,\dots,\theta_n)$ be the complex distance between
$X$ and $Z$.
Denote by $\cal H$ the convex hull of all the vectors in $\R^n$
 obtained from
$\Phi$ by permutations of the coordinates.
Then $\Theta\in\Psi+\cal  H$.}

{\bf 1.6. Riemannian metrics.}
The  $\GL(n,\R)$-invariant
Riemannian metric on $\E_n$
is given by the formula
\begin{equation}
ds^2=\tr\,(dX \cdot X^{-1}\cdot  dX \cdot  X^{-1})
.\end{equation}

{\bf 1.7. Geodesics.}

  {\sc Theorem 1.3.}
{\it Any geodesic in $\E_n$ (or $\PE_n$) has the form
$$
\gamma(t)=
g
\begin{pmatrix}
e^{\phi_1 t}&0&\dots& 0\\
0&e^{\phi_2 t}&\dots& 0\\
\vdots & \vdots& \ddots&\vdots\\
0&0&\dots &e^{\phi_n t}
\end{pmatrix} g^\top
,$$
where $g\in\GL(n,\R)$ and $\phi_1$, \dots, $\phi_n$
are fixed.}

{\sc Corollary 1.4.}
{\it Geodesic distance between $X,Y\in \E_n$
is}
\begin{equation}
\rho(X,Y)=
\Bigl[\sum\nolimits_j \psi_j^2(X,Y)\Bigr]^{1/2}
.\end{equation}

Let normalize the complex distance
in $\PE_n$  by the condition
$\sum\psi_j(X,Y)=0$.
Then the geodesic distance is given by the same formula (1.4).

{\bf 1.8. Cartan subspaces.}
 Cartan subspaces in $\E_n$
are subspaces of the form
$$
L(t_1,\dots,t_n)=
g
\begin{pmatrix}
e^{t_1 }&0&\dots& 0\\
0&e^{t_2 }&\dots& 0\\
\vdots & \vdots& \ddots&\vdots\\
0&0&\dots &e^{t_n }
\end{pmatrix}
 g^\top
,$$
where $g\in\GL(n,\R)$ is fixed, and $t_j$ ranges in $\R$.
Cartan subspaces are totally geodesic submanifolds,
The restriction of the Riemannian metric (1.3)
 to the Cartan subspace $L(t_1,\dots,t_n)$
is $\sum dt_j^2$. In particular, $L(t_1,\dots,t_n)$
is flat, and $t_j$ are flat coordinates.

Let us prove this. Let $\sigma$ ranges in diagonal matrices
whose eigenvalues are $\pm 1$.
Then the maps $X\mapsto \sigma X\sigma^\top$
are involutions, and hence
sets of their fixed points are
totally geodesic submanifolds.
But fixed points for all the such maps are diagonal
matrices. This also implies Theorem 1.3.

  In the language of 1.4, a Cartan subspace consists
of coaxial ellipsoids.

\medskip

{\bf 2. Satake--Furstenberg  boundary}

\medskip

  \addtocounter{sec}{1}
 \setcounter{equation}{0}

\def\1{{\sf Fig. 1.}
\begin{picture}(400,50)(-200,-30)
\put(0,0)
 {\qbezier(-200,0)(-200,10)(0,10)}
{\qbezier(200,0)(200,10)(0,10)}
{\qbezier(-200,0)(-200,-10)(0,-10)}
{\qbezier(200,0)(200,-10)(0,-10)}
{\qbezier(-5,0)(-5,10)(0,10)}
{\qbezier(-5,0)(-5,-10)(0,-10)}
{\qbezier[10](5,0)(5,10)(0,10)}
{\qbezier[10](5,0)(5,-10)(0,-10)}
\multiput(0,0)(0,3){4}{\line(0,2){2}}
\put(0,10){\vector(0,1){10}}
\multiput(0,0)(0,-3){4}{\line(0,-2){2}}
\put(0,-10){\line(0,-1){10}}
\put(5,20){$\scriptstyle Oz$}
\put(4,-20){$\scriptstyle-1$}
\put(0,-10){\circle*{2}}
\multiput(0,0)(3,0){70}{\line(2,0){2}}
\put(200,0){\vector(1,0){20}}
\put(220,3){$\scriptstyle Oy$}
\put(201,3){$\scriptstyle j$}
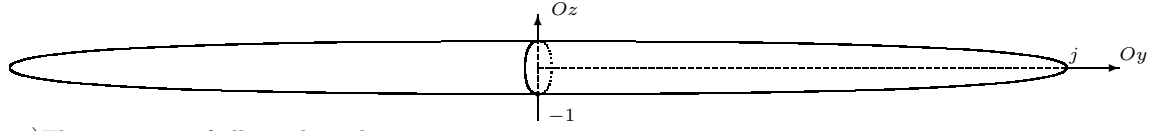
\put(-200,-30){ a)\small The sequence of ellipsoids
with semiaxes $1$, $1$, $j$.}
\end{picture}

\nopagebreak

\begin{picture}(400,50)(-200,-20)
\put(-200,0){\qbezier(-5,0)(-5,10)(0,10)}
\put(-200,0){\qbezier(-5,0)(-5,-10)(0,-10)}
\put(-200,0){\qbezier[10](5,0)(5,10)(0,10)}
\put(-200,0){\qbezier[10](5,0)(5,-10)(0,-10)}
\put(200,0){\qbezier(-5,0)(-5,10)(0,10)}
\put(200,0){\qbezier(-5,0)(-5,-10)(0,-10)}
\put(200,0){\qbezier(5,0)(5,10)(0,10)}
\put(200,0){\qbezier(5,0)(5,-10)(0,-10)}
\put(0,0){\qbezier(-5,0)(-5,10)(0,10)}
\put(0,0){\qbezier(-5,0)(-5,-10)(0,-10)}
\put(0,0){\qbezier[10](5,0)(5,10)(0,10)}
\put(0,0){\qbezier[10](5,0)(5,-10)(0,-10)}
\put(-200,10){\line(1,0){400}}
\put(-200,-10){\line(1,0){400}}
\multiput(0,0)(3,0){70}{\line(2,0){2}}
\put(195,0){\vector(1,0){25}}
\put(220,3){$\scriptstyle Oy$}
\multiput(0,0)(0,3){4}{\line(0,2){2}}
\put(0,10){\vector(0,1){10}}
\multiput(0,0)(0,-3){4}{\line(0,-2){2}}
\put(0,-10){\line(0,-1){10}}
\put(0,-10){\circle*{2}}
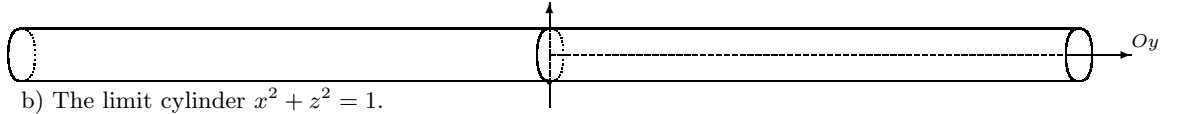
\put(-200,-20){\small b) The limit cylinder $x^2+z^2=1$.}
\end{picture}

\nopagebreak

\begin{picture}(400,50)(-200,-20)
{\linethickness{2pt}
\put(-200,0){\vector(1,0){420}}}
\put(220,3){$\scriptstyle Oy$}
\put(0,0){\qbezier(-5,0)(-5,10)(0,10)}
\put(0,0){\qbezier(-5,0)(-5,-10)(0,-10)}
\put(0,0){\qbezier[10](5,0)(5,10)(0,10)}
\put(0,0){\qbezier[10](5,0)(5,-10)(0,-10)}
\put(-200,-20)
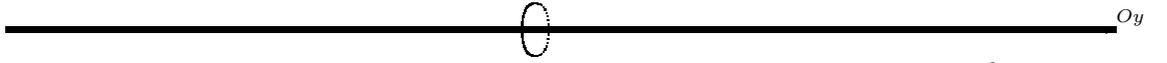
{\small c) The limit in the Semple--Satake space:
the line $Oy$ and the circle in the quotient
space $\R^3/Oy$.}
\end{picture}
}


\def\2{

{\sf Fig. 2.}
\begin{picture}(400,60)(-200,-40)
\put(0,0)
 {\qbezier(-200,0)(-200,20)(0,20)}
{\qbezier(200,0)(200,20)(0,20)}
{\qbezier(-200,0)(-200,-20)(0,-20)}
{\qbezier(200,0)(200,-20)(0,-20)}
{\qbezier(-2,0)(-2,20)(0,20)}
{\qbezier(-2,0)(-2,-20)(0,-20)}
{\qbezier[20](2,0)(2,20)(0,20)}
{\qbezier[20](2,0)(2,-20)(0,-20)}
\multiput(0,0)(0,3){8}{\line(0,2){2}}
\put(0,20){\vector(0,1){10}}
\multiput(0,0)(0,-3){8}{\line(0,-2){2}}
\put(0,-20){\line(0,-1){10}}
\put(5,35){$\scriptstyle Oz$}
\put(4,-27){$\scriptstyle-j$}
\put(0,-20){\circle*{2}}
\multiput(0,0)(3,0){70}{\line(2,0){2}}
\put(200,0){\vector(1,0){20}}
\put(220,3){$\scriptstyle Oy$}
\put(201,3){$\scriptstyle j^2$}
\put(0,0){\vector(-1,-1){10}}
\put(-10,-3){\small 1}
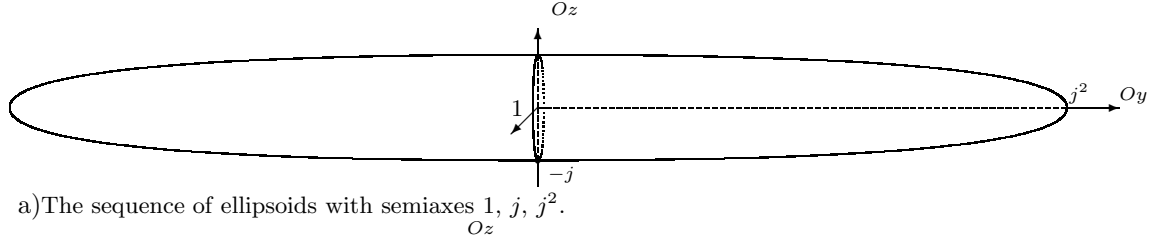
\put(-200,-40){ a)\small The sequence of ellipsoids
with semiaxes $1$, $j$, $j^2$.}
\end{picture}

\nopagebreak

\begin{picture}(400,80)(-200,-40)
\multiput(0,0)(0,3){8}{\line(0,2){2}}
\put(0,20){\vector(0,1){10}}
\multiput(0,0)(0,-3){8}{\line(0,-2){2}}
\put(0,-22){\line(0,-1){10}}
\put(5,35){$\scriptstyle Oz$}
\put(4,-27){$\scriptstyle-j$}
\multiput(0,0)(3,0){70}{\line(2,0){2}}
\put(200,0){\vector(1,0){20}}
\put(220,3){$\scriptstyle Oy$}
\put(0,0){\vector(-1,-1){10}}
\put(-202,-22){\line(1,0){400}}
\put(-202,-22){\line(0,1){40}}
\put(-202,18){\line(1,0){400}}
\put(198,-22){\line(0,1){40}}
\put(-198,22){\line(1,0){400}}
\put(-198,22){\line(0,-1){4}}
\put(202,22){\line(0,-1){40}}
\multiput(-198,-18)(6,0){67}{\line(1,0){3}}
\multiput(-198,18)(0,-6){6}{\line(0,-1){3}}

\put(-200,-40){\small b) The pair of planes $x^2=1$.}

\end{picture}

\nopagebreak

\begin{picture}(400,80)(-200,-40)

{\linethickness{2pt}
\put(-200,0){\vector(1,0){420}}}
\put(220,3){$\scriptstyle Oy$}
\put(0,0){\vector(0,1){40}}
\put(0,0){\line(0,-1){30}}
\put(5,40){$\scriptstyle Oz$}
{\linethickness{1pt}
\put(0,0){\line(0,1){20}}
\put(0,20){\line(1,0){200}}
\put(200,0){\line(0,1){20}}
}
\put(0,0){\vector(-1,-1){15}}
\put(-200,-40)
{\small c) The limit in the Semple--Satake space is
the flag: the line $Oy$
and the plane $yOz$.}
\end{picture}
}

The object described in this section
 is  called the Satake--Furstenberg
compactification
of a symmetric space (see \cite{Sat}).
Its explicit construction given below
arises to Semple \cite{Sem1}, \cite{Sem3}
and Alguneid \cite{Alg}.

{\bf 2.1. Semple--Satake space.}
A point of the Semple--Satake space $\S_n$ is the following
collection (i)--(iii) of the data.

\smallskip

(i) A subset
\begin{equation}
I=\{i_1,i_2,\dots,i_p\}\subset
\{1,2,\dots, n-1\}
,
\end{equation}
where $p=0,1,\dots,n-1$.  It is convenient to assume
$i_0=0$, $i_{p+1}=n$.

\smallskip

(ii) A flag
\begin{equation}
\R^n=W_0\supset W_1\supset W_2\supset
\dots \supset W_p\supset W_{p+1}=0
,\end{equation}
 where  for each $k=1,\dots,p$
 \begin{equation}
\codim W_k=i_k
.\end{equation}

\smallskip

(iii) A collection
$$Q_j; \qquad j=1,\dots,p+1$$
of  positive definite quadratic forms
 on $W_{j-1}/W_j$ defined up to  scalar factors%
\footnote{equivalently, we have an  ellipsoid
defined up to a homothety in each
subquotient $W_{j-1}/W_j$.}.

\medskip

We denote the piece of $\S_n$
corresponding to the collection (2.1) by
$$\S_n(I)=\S_n(i_1,i_2,\dots,i_p).$$
Thus,
$$
\S_n=\bigcup_{I\subset \{1,2,\dots,n-1\}}
     \S_n(I)
.$$
We also denote by $\F(I)$ the set of all the flags (2.2)
satisfying (2.3).

\smallskip

{\sc Remark.} The set $\S_n(\emptyset)$
is $\PE_n$.
A set $\F(I)$ is a fiber bundle whose
base is the space $\F(I)$ and fibers are
\begin{equation}
\prod\nolimits_{k=1}^{p+1}\PE_{i_{k}-i_{k-1}}
.\end{equation}
Each expert in semisimple groups
can easily translate this form of definition into
the root language.

\smallskip

{\sc Remark.} In particular, for $I=\{1,2,\dots,n-1\}$,
 the fibers are points, and hence
$\S_n(I)=\F(I)$
is the space of complete  flags.

A simple calculation shows that
 $$\dim \S_n(i_1,\dots i_p)=n^2-1-p=
\dim\PE_n-p
.$$

Now  we will define the topology
of a compact metrizable space
on $\S_n$, this topology satisfies
the property:
$$
\text{
Closure of} \quad
\S_n(I)=\bigcup\nolimits_{J\supset I} \S_n(J)
.$$
In particular, the closure
of $\S_n(\emptyset)= \PE_n$ is
the whole space
$\S_n$.

{\bf 2.2. Inductive definition of  convergence in $\S_n$.}
Assume that the convergence is defined
in all the spaces $\S_m$ for $m<n$.
Consider a sequence of $n$-dimensional ellipsoids
$$Q^{(j)}:\quad \tfrac12 \sum\nolimits_{k,l} x_{kl}^{(j)} v_k v_l=1$$
defined up to a homothety.
We can assume that the shortest semiaxis of each ellipsoid
$Q^{(j)}$ equals 1.

The first necessary condition of convergence is the
convergence of ellipsoids (normalized in this way)
in the human sense, i.e.,
 the convergence of
the corresponding matrices $X^{(j)}$.
Denote the limit by $Y$.

If $Y$ is a nondegenerate matrix, then $Y$
is the a limit
in the sense of the space $\S_n$
(and $Y\in\PE_n$).

Assume that the matrix $Y$ is degenerate.
Consider the kernel $\ker Y$, denote its dimension by $m$;
geometrically
$\ker Y$ is the directing subspace of the cylinder
$\tfrac12\sum y_{kl}v_kv_l=1$.
In particular, we obtain the ellipsoid
(defined up to a homothety) in the quotient
space $\R^n/ \ker Y$
and the sequence of ellipsoids
$
Q^{(j)}\cap \ker Y
$
in the subspace $\ker Y$.
Now, the sufficient and necessary condition
of the convergence
is the convergence of the sequence   $Q^{(j)}\cap \ker Y $
of ellipsoids
in the sense of the space $\S_m$ (where $m=\dim\ker Y)$.

{\bf 2.3. Examples.}
Let $n=3$, i.e., we have a sequence of ellipsoids
in $\R^3$.

{\sc Example 1,}
see Fig. 1.
Consider the sequence of ellipsoids
$$
Q^{(j)}:\,\,
x^2+{y^2}/{j^2}+z^2=1
.$$

If $j\to\infty$, this family of surfaces
converges  to the cylinder
$x^2+z^2=1$. Its axis is $Oy$.
Thus we obtain the flag $\R^3\supset Oy\supset 0$
and the  circle $x^2+z^2=1$ in the quotient
$\R^3/Oy$.

  \begin{figure}
\1
\end{figure}

{\sc Example 2.}
Consider the sequence
of ellipsoids
$$Q^{(j)}:\,\,
x^2+{y^2}/{j^4}+ {z^2}/{j^2}=1
.$$
The limit is the pair of planes
$x^2=1$.
The directing plane of $x^2=1$  is
$yOz$.
The section of $Q^{(j)}$ by the plane
$x=0$ is the ellipse
$
{y^2}/{j^4}+{z^2}/{j^2}=1$.
We consider quadrics up to a homothety,
 hence we can replace our equation by
$
{y^2}/{j^2}+z^2=1
$
This sequence of ellipses converges to the pair
of lines
$z^2=1$.
Its directing line is
$z=0$.
Finally, we obtain the flag
$
\R^3\supset yOz\supset Oy
$ as the limit of the sequence of ellipsoids.

\begin{figure}
\2
\end{figure}

{\bf 2.4. Noninductive definition of convergence in $\S_n$}
(it is used only in 3.9).
Consider a sequence of positive definite matrices
$X^{(j)}\in\PE_n$.
Denote their eigenvalues   by
\begin{equation}
\lambda_1^{(j)}\ge
\lambda_2^{(j)}\ge
\dots\ge
\lambda_n^{(j)}
.\end{equation}

Now we present the conditions for
the sequence $X^{(j)}$ be convergent.

\smallskip

{\bf Condition A.}
There exists a separation of the set of
eigenvalues (2.5) into the 'packets'
$$ \Bigl(\alpha_1^{(j)},\dots,
    \alpha_s^{(j)}\Bigr),\,
 \Bigl(\beta_1^{(j)},\dots,
    \beta_t^{(j)}\Bigr),\,
\Bigl(\gamma_1^{(j)},\dots,
    \gamma_r^{(j)}\Bigr),\,
\Bigl(\delta_1^{(j)},\dots,
  \delta_p^{(j)}\Bigr),\dots
$$
such that $s$, $t$, $r$, etc. are independent on
$j$ and the following list of the conditions
$1^\circ-3^\circ$ is satisfied.

\smallskip

 $1^\circ$. For sufficiently large values of $j$,
$$\alpha_1^{(j)}\ge\dots\ge
    \alpha_s^{(j)}\ge
\beta_1^{(j)}\ge\dots\ge
    \beta_t^{(j)}\ge\gamma_1^{(j)}\ge\dots
$$

 $2^\circ$. For the eigenvalues from one packet
we have
\begin{gather*}
\forall k,m \qquad
\lim\limits_{j\to\infty}
{\alpha_k^{(j)}}/{\alpha_m^{(j)}}\qquad
\text{is finite and nonzero}\\
\forall u,v \qquad
\lim\limits_{j\to\infty}
{\beta_u^{(j)}}/{\beta_v^{(j)}}\qquad
\text{is finite and nonzero}
\end{gather*}
etc.

$3^\circ$. For the eigenvalues from different packets
we have
\begin{gather*}
\forall k,u \qquad
\lim\limits_{j\to\infty}
{\alpha_k^{(j)}}/{\beta_u^{(j)}}=\infty,
\qquad\quad
\forall u,v \qquad
\lim\limits_{j\to\infty}
{\beta_u^{(j)}}/{\gamma_v^{(j)}}=\infty\\
\forall v,w \qquad
\lim\limits_{j\to\infty}
{\gamma_v^{(j)}}/{\delta_w^{(j)}}=\infty,
\qquad \text{etc.}
\end{gather*}

{\bf Condition B.}
Denote by the $V_{\alpha}^{(j)}\subset \R^n$ the subspace
spanned by the eigenvectors corresponding
 to the eigenvalues
$\alpha_1^{(j)},\dots,
    \alpha_s^{(j)}$ of $X^{(j)}$; in the same way we define
$V_\beta^{(j)}$, $V_\gamma^{(j)}$, etc%
\footnote{For  sufficiently large
values of $j$ we have
$\alpha_s^{(j)}>\beta_1^{(j)}$,
and hence $V_\alpha$ is well-defined.}.
Consider the subspaces
\begin{gather*}
W_1^{(j)}=V_\beta^{(j)}\oplus
V_\gamma^{(j)}\oplus
V_\delta^{(j)}\oplus\dots\\
W_2^{(j)}=\hphantom{V_\beta^{(j)}\oplus}
V_\gamma^{(j)}\oplus
V_\delta^{(j)}\oplus\dots\\
W_3^{(j)}=\hphantom{Y_\beta^{(j)}\oplus
V_\gamma^{(j)}\oplus}
V_\delta^{(j)}\oplus\dots
\end{gather*}
etc.
Our requirement is

{\it for each $q$ the sequence
 of subspaces $W_q^{(j)}$ converges to some $W_q$.}

Thus, we obtain the flag
$\R^n=W_0\supset W_1\supset W_2\supset \dots $
and our next purpose is to obtain
a quadratic form in each subquotient
$W_q/W_{q+1}$.

{\bf Condition C.} Denote by $R_q^{(j)}$
the restriction of the bilinear form
$$
Q_X^j(v,w)=\tfrac12 \sum x_{kl}^{(j)} v_k w_l
$$
to the subspace $W_q^{(j)}$.
These forms are defined up to  scalar factors,
we fix these factors
from the condition: the shortest
semiaxis of $R_q^{(j)}$ is 1.

Our last requirement is:

{\it  for each $q$ the sequence $R_q^{(j)}$ of
bilinear forms converges as $j\to\infty$.}

We must say this more carefully, since the forms
$R_q^{(j)}$
are defined on different subspaces.
For this, consider a sequence of orthogonal operators
$h^{(j)}\in\SO(n)$ such that $h^{(j)}$
converges to $E$ and $h^{(j)} W_q^{(j)}=W_q$.
Thus we identify $W_q^{(j)}$ and $W_q$.
After this we can tell about convergence of
bilinear forms on $W_q$.

 Denote by $R^\times_q$ the limit of the sequence
$R_q^{(j)}$. Evidently,
$
R^\times_q(v,w)=0$
for
$v\in W_q,\,\, w\in W_{q+1}
$.
Hence, we obtain the well-defined
bilinear form $R^\square_q$ on the quotient
space $W_q/W_{q+1}$.

Thus, the convergence of a sequence in $\PE_n$
to a point of $\S_n$ is defined.

{\bf 2.5. Convergence on the boundary.}
As we have seen, each set $\S_n(I)$
is a bundle over $\F(I)$ with fibers (2.4).
Let us compactify each factor
$\PE_{i_k-i_{k+1}}$ as $\S_{i_k-i_{k+1}}$.
 Thus we obtain a compactification
$\overline{\S_n(I)}$ of $\S_n(I)$.

But we have
the obvious embedding $\overline{\S_n(I)}\to\S_n$.
This remark also defines the convergence on the boundary.

{\bf 2.6. Result.}
{\sc Theorem 2.1.}
{\it $\S_n$ is a compact metrizable topological space.}

\medskip

{\bf 3. Pencils of geodesics}

\medskip

  \addtocounter{sec}{1}
 \setcounter{equation}{0}

\nopagebreak

The term 'geodesic' here and below
means a {\it directed geodesic}.

Following Karpelevich, in 3.2--3.6 we define
 and describe explicitly 3 types of pencils of
geodesics in $\E_n$ and $\PE_n$.
In 3.7-3.9 we define limit points of  pencils at infinity.
In Section 5,
we also need  in description of
pencils in products of the type
$\PE_{k_1}\times\dots\times \PE_{k_m}$.
The necessary modification
of the constructions of pencils is given in
3.10-3.12.

\smallskip

{\bf A. Definitions and canonical forms of pencils}

\smallskip

\nopagebreak

{\bf 3.1. Velocities of geodesics.}
Consider a directed geodesic
\begin{equation}
\mu(t)=g\begin{pmatrix}
e^{\phi_1 t} &0&\dots\\
0&e^{\phi_2 t}&\dots\\
\vdots&\vdots &\ddots
\end{pmatrix}
g^\top
\end{equation}
with
\begin{equation}
\phi_1\ge\phi_2\ge\dots\ge\phi_n
.\end{equation}
Its {\it velocity} in $\E_n$ is the
collection of numbers
(3.2)
defined up to a joint positive factor;
this freedom corresponds to the substitution
$t=a t'$ to (3.1).

The velocity of $\mu(t)$ in $\PE_n$
is the vector (3.2) defined up to
transformations
\begin{equation}
(\phi_1,\phi_2,\dots)\mapsto
(a\phi_1+b,a\phi_1+b,\dots)
.\end{equation}

We need in an overfilled system of notation
for velocity vectors.

Fix positive integers
$\alpha_1$, \dots, $\alpha_m$
such that $\sum\alpha_j=n$.
Fix real numbers
\begin{equation}
\psi_1>\psi_2>\dots>\psi_m
.\end{equation}
For such data we compose the velocity vector
\begin{equation}
(\underbrace
{\psi_1,\dots,\psi_1}_{ \alpha_1\,\, \text{times}},
\underbrace
{\psi_2,\dots,\psi_2}_{ \alpha_2\,\, \text{times}},
\underbrace
{\psi_3,\dots,\psi_3}_{ \alpha_3\,\, \text{times}},
\dots)
.\end{equation}
In this notation, the geodesic (3.1) can be written in the form
\begin{equation}
\mu(t)=
g\begin{pmatrix}
e^{\psi_1 t}E_{\alpha_1} &0&\dots\\
0&e^{\psi_2 t}E_{\alpha_2}&\dots\\
\vdots&\vdots &\ddots
\end{pmatrix}g^\top
,\end{equation}
where $E_\alpha$ denotes the unit
$\alpha\times\alpha$ matrix.

We also define the subset
\begin{equation}
I=\{i_0,i_1,\dots,i_m\}\subset \{0,1,2,\dots,n\}
\end{equation}
by
\begin{align}
&i_0=0,\nonumber\\
&i_k=\alpha_1+\dots+\alpha_k \qquad \text{for $k=1,\dots,m-1$}\\
&i_m=\alpha_1+\dots+\alpha_m=n    \nonumber
\end{align}


For any subset $I=\{0,i_1,\dots,i_{m-1},n\}\subset \{0,1,2,\dots,n\}$
 denote by $\Delta(I)$ the simplex consisting of
 collections
(3.4) defined up to a positive factor;
 by $\Delta_\circ(I)$
 denote the set of all collections (3.4) defined up
to the equivalence (3.3).

\smallskip

We defined the velocity
of a geodesic in the terms of its canonical form.
Let us define it in terms of complex distance.

{\sc Lemma 3.1.} {\it
Denote by
$
e^{\tau_1(t)}\ge e^{\tau_2(t)}\ge\dots
$
the eigenvalues of the matrix $\mu(t)$ given by (3.1).
 For all $j$, we have}
$\lim\limits_{t\to+\infty} {\tau_j(t)}/{t}=\phi_j.$

{\sc Corollary 3.2.}
{\it Fix $A\in\E_n$. Denote by $\sigma_1(t)\ge\sigma_2(t)\ge\dots$
the complex distance between $A$ and $\mu(t)$.
Then for all $j$ we have}
$\lim\limits_{t\to+\infty} {\sigma_j(t)}/{t}=\phi_j$.

These statements follow from Theorem 1.2.

{\bf 3.2. Null
pencils and finite pencils.}
Consider a {\it directed} geodesic $\mu(t)$
(in $\E_n$ or $\PE_n$).
The corresponding {\it null pencil}
$\Pi^\nul_\mu$ is the set of all
the geodesics $\nu(t)$ such that
$$
\lim\limits_{t\to +\infty}\dist(\mu(t),\nu)=0,$$
where the distance between a point $X$ and a geodesic $\nu$
is
$$
\dist(X,\nu):=\min\limits_{s\in\R}
\rho(X,\nu(s))
$$

We also define the {\it finite pencil} $\Pi^\fin_\mu$
as the set of all the geodesics
$\nu(t)$ such that there exists
a finite limit
$$\lim_{t\to+\infty} \dist(\mu(t),\nu)$$

{\bf 3.3. Canonical forms of null pencils
and finite pencils.}
Consider a geodesic $\gamma$
given by
\begin{equation}
\gamma(t)=
\begin{pmatrix}
e^{\psi_1 t}E_{\alpha_1} &0&\dots\\
0&e^{\psi_2 t}E_{\alpha_2}&\dots\\
\vdots&\vdots &\ddots
\end{pmatrix}
,\end{equation}
numbers $\psi_j$ satisfy (3.4).

{\sc Theorem 3.3.}
 a) {\it The finite pencil
$\Pi^\fin_\gamma$ consists of all the geodesics
that can be represented in the form
\begin{equation}
h\gamma(t)h^\top
,\end{equation}
where $h$ ranges in
 $(\alpha_1+\alpha_2+\dots)\times
(\alpha_1+\alpha_2+\dots)$ block matrices of the shape}
\begin{equation}
h=
\begin{pmatrix}
H_{11}&0&0&\dots\\
H_{12}&H_{22}&0&\dots\\
H_{13}&H_{23}&H_{33}&\dots\\
\vdots&\vdots&\vdots&\ddots
\end{pmatrix}
.\end{equation}

 b) {\it The null pencil
$\Pi^\nul_\gamma$ consists of all the geodesics
having the form (3.10),
where $h$ ranges in
 $(\alpha_1+\alpha_2+\dots)\times
(\alpha_1+\alpha_2+\dots)$ block matrices of the shape}
\begin{equation}
h=
\begin{pmatrix}
E_{\alpha_1}&0&0&\dots\\
H_{12}& E_{\alpha_2} &0&\dots\\
H_{13}&H_{23}& E_{\alpha_3} &\dots\\
\vdots&\vdots&\vdots&\ddots
\end{pmatrix}
.\end{equation}

{\sc Proof.} Let $\mu\in\Pi^\fin_\gamma$.
Theorem 1.2 implies coincidence of the
velocities of $\mu$, $\gamma$,
 hence $\mu(t)=g\gamma(t)g^\top$ for some
$g\in\GL(n,\R)$. By the same theorem, for large $|t-s|$
the points $\mu(s)$ and $\nu(t)$ are far.

  Let $\lambda_k(t)=e^{\psi(t)}$ be the solutions of the equation
$\det(\lambda\gamma(t)-g\gamma(t)g^\top)=0$, i.e., $\psi_k(t)$
is the complex distance between $\gamma(t)$ and $\mu(t)$.
Equivalently, $\lambda_k(t)$ are the eigenvalues
of
$\gamma(t)^{-1/2} g \gamma(t) g^\top \gamma(t)^{-1/2}$.
Thus, $\lambda_k(t)^{1/2}$ are the singular values of
$Z(t):=\gamma(t)^{-1/2} g\gamma(t)^{1/2}$. But the numbers
$\psi_k(t)$ are bounded, hence  matrix elements
of $Z(t)$ are bounded. Therefore, $g$
is triangular.    \hfill $\square$

All geodesics lying in a given pencil
have the same velocity. Thus the {\it velocity
of a pencil} is well defined.

{\bf 3.4. Decomposition of finite pencils
 into null pencils.}
Denote by $P_\gamma$ the group of all the
matrices (3.11), and by $N_\gamma$ the group
of all the matrices
(3.12). Evidently, $N_\gamma$ is
 a normal subgroup in $P_\gamma$.

The group $P_\gamma$ acts on $\Pi_\gamma^\fin$
by the transformations
$h: \mu(t)\mapsto h\mu(t)h^\top$
and the subgroup $N_\gamma$
transfers the null pencil   $\Pi^\nul_\gamma$
 to itself. Moreover, $N_\gamma$ transfers
each null subpencil of $\Pi^{fin}_\gamma$
to itself.

Let $\mu_1$, $\mu_2\in \Pi^\fin_\gamma$.
We say that
$\mu_1\sim\mu_2$
if they lie in one null pencil.
We denote by $\widetilde\Pi_\gamma$
the quotient of  $ \Pi^\fin_\gamma$
by
this equivalence relation.
The group
$$
P_\gamma/N_\gamma\simeq
\GL(\alpha_1,\R)\times\GL(\alpha_2,\R)\times\dots
$$
acts on $\widetilde\Pi_\gamma$ in a natural way,
and we obtain
\begin{equation}
\widetilde\Pi_\gamma\simeq
\GL(\alpha_1,\R)/\O(\alpha_1)
\times\GL(\alpha_2,\R)/\O(\alpha_2)\times\dots
=\E_{\alpha_1}\times \E_{\alpha_2}\times\dots
\end{equation}

{\sc Remark.}
(\cite{Kar})
For $\nu_1$, $\nu_2\in\Pi_\gamma^\fin$
we can define the distance at infinity
$$
\dist(\nu_1,\nu_2)=
\lim\limits_{t\to\infty}
\inf\limits_{s\in\R}
\rho\bigl[\nu_1(t),\nu_2(s)\bigr]
.$$
Then this distance is the geodesic distance
in the symmetric space (3.13).

{\bf 3.5. Solvable pencils.}
Finite  pencils and null pencils have sense
for any space of nonpositive curvature.
For symmetric spaces there exists
a natural
intermediate  equivalence of
geodesics.

In the previous subsection, we constructed the map
$$
\Pi_\gamma^\fin\to
\prod\nolimits_j \GL(\alpha_j,\R)/\O(\alpha_j)
.$$

The symmetric
space $\prod \GL(\alpha_j,\R)/\O(\alpha_j)\simeq \prod\E_{\alpha_j}$
 is not semisimple. Consider
its natural projection
to the semisimple space
$
\prod\nolimits_j \PGL(\alpha_j,\R)/\PO(\alpha_j)
\simeq \prod\nolimits_j\PE_{\alpha_j}
$.
Thus, we obtain the canonical map
\begin{equation}
\Pi^\fin_\gamma\to\prod\nolimits_j \PGL(\alpha_j,\R)/\PO(\alpha_j)
\simeq\prod\nolimits_j\PE_{\alpha_j}
.\end{equation}

We say that $\nu_1$, $\nu_2$ are elements of one
{\it solvable pencil} if their images
under this map coincide.

For any geodesic $\gamma$, we denote by $\Pi_\gamma^\sol$
the corresponding solvable pencil. Obviously, we have
$$
\Pi_\gamma^\fin\supset
\Pi_\gamma^\sol
\supset
\Pi_\gamma^\nul
.$$

{\bf 3.6. Canonical forms of solvable pencils.}

{\sc Proposition 3.4.}
{\it
Let a geodesic $\gamma$ has the form (3.9). Then the corresponding
solvable pencil $\Pi_\gamma^\sol$ consists of all the geodesics
 having the form $h\gamma(t)h^\top$,
where $h$ is an
 $(\alpha_1+\alpha_2+\dots)\times
(\alpha_1+\alpha_2+\dots)$ block matrix of the shape
\begin{equation}
h=
\begin{pmatrix}
\tau_1 E_{\alpha_1}&0&0&\dots\\
H_{12}& \tau_2 E_{\alpha_2} &0&\dots\\
H_{13}&H_{23}& \tau_3 E_{\alpha_3} &\dots\\
\vdots&\vdots&\vdots&\ddots
\end{pmatrix}
,\qquad\qquad \tau_j\in\R
.\end{equation}
}

The set of solvable pencils in a given finite pencil
(3.10)--(3.11)
is parametrized by the collection of diagonal blocks
$H_{jj} H_{jj}^\top$
of the matrix (3.10), these blocks are defined up
to a multiplication by positive scalars.

\smallskip

{\bf B. Boundary data for pencils}

\smallskip

\nopagebreak

{\bf 3.7. Boundary data for solvable pencils.}

\nopagebreak

{\sc Lemma 3.5.} a) {\it Each geodesic $\mu(t)$
has a limit in $\S_n$ as $t\to\infty$.}

b) {\it If $\mu_1$, $\mu_2$ lie in one solvable
pencil, then their limits in $\S_n$ coincide.}

In fact, evaluation of the limit of a geodesic in $\S_n$
  is reduced
to  evaluation of limit of a family of coaxial
ellipsoids.
For instance, let us describe explicitly the
 limit of a geodesic $\gamma$
given by (3.9)
. Denote by $e_j$ the standard basis
in $\R^n$.
Denote by $W_k$ the subspace
$
W_k=\oplus_{p>i_k} \R e_p
$
Thus we obtain the flag
$\R^n=W_0\supset W_1 \supset \dots$.
The ellipsoids in the quotients $W_{k-1}/W_k$
are spheres.

Now let $\mu(t)$ be another geodesic of the same
finite pencil, i.e.,  $\mu(t)=h\gamma(t)h^\top$
with $h$ given by (3.11). Then the limit flag is the same,
and
the quadratic forms
in the quotients $W_{k-1}/W_k$
are
$Q_k=H_{kk}H_{kk}^\top$,
they are defined up to a multiplication by a scalar factor.

Now the following statement becomes obvious.

{\sc Theorem 3.6.}
{\it A solvable pencil is uniquely determined by
its velocity and its limit in the space
 $\S_n$.}

Consider a geodesic $\mu(t)$,
whose velocity
 is contained in the simplex $\Delta(I)$, see 3.1.
As we have explained above,
$\lim \mu(t)$ in the sense of $\S_n$ belongs $\S_n(I)$;
 the strata $\S_n(I)$ were defined
 in 2.1.

{\sc Corollary 3.7.}
{\it Fix a set $I$. Denote by $\K(I)$
the space of all the solvable pencils in $\PE_n$, whose
velocities have the form (3.5).
Then
$$
\K(I)
\simeq
\Delta_\circ(I)\times
\S_n(I)
$$
and}
\begin{equation}
\dim \K(I)=n^2-2=\dim\PE_n-1
\end{equation}


{\bf 3.8. Boundary data for finite pencils.}
As we have seen, each geodesic $\mu(t)$
has a limit in $\S_n$. In particular,
we have a canonically defined flag in $\R^n$,
we call it by the {\it limit flag}.

{\sc Proposition 3.8.}
{\it Two geodesics lie in one finite pencil
iff their limit flags and their velocities
coincide.}


%


Thus,
for any finite pencil we associate the following collection
of data.

1) A set
$I=\{0,i_1,\dots, i_{k-1},n\}\subset\{0,1,2,\dots,n\}$, where $k> 1$.

2) A point of the simplex $\Delta(I)$ (or $\Delta_\circ(I)$
for $\PE_n$).

3) A flag lying in $\F(I)$ (see Subsection 2.1).

\smallskip

{\bf 3.9. Boundary data for null pencils in $\E_n$.}
Obviously, in this case we must remember more than
in the case of solvable pencils.

Consider a geodesic $\mu\subset \E_n$,
whose velocity has the form (3.5).
Let us fix a parameter $t$ on $\mu$.
This means that we fix an origin $X_0\in \mu$
and we fix the
velocity vector $(\psi_1,\psi_2,\dots)$ literally
(without any equivalence).

 Preserving the notation of 2.1,
we  change the construction
of 2.4 in one place.

Thus, we have the family of quadratic forms
$Q(t)$ on $\R^n$ corresponding to points of the geodesic $\mu
(t)$.
After obtaining the one-parametric family of flags
$
\R^n=W_0(t)\supset W_1(t)\supset W_2(t)\supset\dots,
$
we consider the restriction
$R_p(t)$ of the form $Q(t)$ to the subspace
$W_{p-1}(t)$. Then we
consider
$$
R^\times_p:=\lim\limits_{t\to +\infty}
e^{-\psi_p t} R_p(t)
$$
This is a well defined  nondegenerate symmetric bilinear form
on $W_{p-1}$. The subspace
$W_p$ is the kernel of this form,
and finally we obtain a nondegenerate symmetric bilinear form
$R^\square_p$ on each subquotient
$W_{p-1}/W_p$.

We emphasis, that in Section 2 the forms $R^\square_p$
on subquotients were defined up to positive factors.
Now they are defined literally.

\smallskip

But we started from a geodesic with
a fixed parametrization.

If we  multiply the velocity $(\psi_1,\psi_2,\dots)$
by a scalar and leave the origin,
then our limit data (the flag $W_p$ and the forms
$R^\square_p$) do not  change.

If we move the origin $X_0$ along the geodesic,
then the collection $R^\square_p$ changes
in the following way
\begin{equation}
(R_1^\square,R_2^\square,\dots)
\rightsquigarrow
(e^{\psi_1 s}R_1^\square,e^{\psi_2 s}R_2^\square,\dots)
\qquad \text{for some $s\in \R$}.
\end{equation}

Thus for any null pencil, we associate
the following boundary data

1) A  set
$I=\{0,i_1,\dots, i_{m-1},n\}\subset\{0,1,2,\dots,n\}$,
where $m>1$.

2) A point of the simplex $\Delta(I)$.

3) A flag
$W_1\supset W_2\supset\dots \supset
W_{m-1}$
lying in $\F(I)$

4) The family
$(R_1^\square,R_2^\square,\dots,R_{k}^\square)$
of positive definite quadratic forms
on subquotients $W_{p-1}/W_p$
defined up to the equivalence (3.17).

{\sc Theorem 3.9.}
{\it The space of all null pencils is
in one to one correspondence
with the collections of data 1-4.}

\smallskip

{\bf C. Pencils in products of symmetric spaces}

\smallskip

{\bf 3.10. Abstract definition of  solvable pencils.}
 Let us give a definition of
a solvable pencil
 in a semisimple  Riemannian symmetric spaces $G/K$;
$K$ is a maximal compact subgroup in $G$.
The only case interesting for us is $G/K=\prod\PE_{m_j}$.

First, consider the group $P_\gamma$ of all the isometries
of $G/K$ mapping
a finite pencil $\Pi_\gamma^\fin$ to itself.
It is a parabolic
 subgroup (\cite{Kar}), and
it acts transitively on $\Pi_\gamma^\fin$.
Denote by $L_\gamma$ the Levi factor in $P_\gamma$
(maximal reductive subgroup)
Denote by $K_\gamma$
the maximal compact subgroup
in $L_\gamma$.

As in 3.5, the space of all null-pencils
in $\Pi^\fin_\gamma$ can be identified with
$L_\gamma/K_\gamma$.
It is a reductive symmetric space,
hence it is
a product
$S_\gamma\times L_\gamma$, where $S_\gamma$ is a semisimple symmetric space
and $L_\gamma$ is an Euclidean space.
We say that two elements of
 $\Pi^\fin_\gamma$ lie in one solvable pencil
if their images under the map
$
\Pi^\fin_\gamma \to L_\gamma/K_\gamma \to S_\gamma
$
coincide.

{\bf 3.11. Vel-geodesics.}
Consider a
 Riemannian manifold $M$.
 We say, that a parameter
$t$ on a geodesic $\gamma$ is semi-natural,
if the Riemannian length of the tangent vector
$\|\gamma'(t)\|=c$
is independent on $t$.
We say, that a {\it vel-geodesic} is a directed geodesic with a fixed
constant $c$.

{\bf 3.12. Pencils in
products of spaces $\E_{m_j}$}.
Consider a
vel-geodesic $\gamma(t)$ in the space
\begin{equation}
\prod\nolimits_{\tau=1}^\beta\PE_{m_\tau}.
\end{equation}
Let $\gamma_\tau(t)$ be   projections
of $\gamma(t)$ to $\PE_{m_\tau}$.


We define a solvable vel-pencil of geodesics
as a solvable pencil with fixed
constant $c$ as in 3.11.

{\sc Lemma 3.10.}
{\it Geodesics $\gamma$, $\mu$ lie in one solvable
 vel-pencils, iff
for each $\tau$  the geodesics $\gamma_\tau$, $\mu_\tau$
lie in one solvable vel-pencil.}

In particular, the space of all solvable pencils
in our space
is parametrized by the following collections of data

A. Family $I_1$, \dots, $I_\beta$
of subsets
$$
I_\tau=
\{0,i_1^{(\tau)},i_2^{(\tau)},\dots, i_{j_\tau}^{(\tau)}, m_\tau \}
\subset
    \{0,1,2,\dots,m_\tau-1, m_\tau\}
$$

B. A collection  of numbers
\begin{equation}
\psi^{(1)}_1>\dots>\psi^{(1)}_{j_1+1}; \qquad
\psi^{(2)}_1>\dots>\psi^{(2)}_{j_2+1};\qquad \dots
.\end{equation}
defined up to the equivalence
\begin{multline}
\Bigl\{
\bigl[\psi^{(1)}_1,\dots,\psi^{(1)}_{j_1+1}\bigr],
\bigl[\psi^{(2)}_1,\dots,\psi^{(2)}_{j_2+1}\bigr],
\bigl[\psi^{(3)}_1,\dots,\psi^{(3)}_{j_3+1}\bigr],
\dots\Bigr\}
\sim
\\
\sim
\Bigl\{
\bigl[a\psi^{(1)}_1+b_1,\dots,a\psi^{(1)}_{j_1+1}+b_1\bigr],
\bigl[a\psi^{(2)}_1+b_2,\dots,a\psi^{(2)}_{j_2+1}+b_2\bigr],
\\
\bigl[a\psi^{(3)}_1+b_3,\dots,a\psi^{(3)}_{j_3+1}+b_3\bigr],
\dots\Bigr\}
.\end{multline}

C. A point of
$\S_{m_1}(I_1)\times \dots \times\S_{m_\beta}(I_\beta)$.

 \smallskip

We denote by $\Delta_\circ(I_1,\dots, I_\tau)$
the set of all the collections (3.19)
defined up to the equivalence (3.20).

\medskip

{\bf 4. Finite pencils. Matrix sky and its tilling}

\medskip

  \addtocounter{sec}{1}
 \setcounter{equation}{0}

Now we want to construct an ideal boundary
of $\E_n$ or $\PE_n$ as a set of limits
of pencils of geodesics. All the three types of pencils are
available for this purpose, but the final results
in these three cases are essentially distinct.

In this section, there is no difference between
$\E_n$ and $\PE_n$. For definiteness, we discuss $\E_n$.

{\bf 4.1. Sphere at infinity.}

{\sc Proposition 4.1.}
{\it For each point $X\in\E_n$
and each finite pencil $\Pi^\fin_\gamma$,
there exists a unique geodesic $\mu\in\Pi^\fin_\gamma$
passing the point $X$.}

{\sc Proof.}  Consider the geodesic $\gamma(t)$
given by (3.9). Denote by $P_\gamma$
the group of all the matrices (3.11).
 First, each positive matrix can be represented
in the form $hh^\top$, where $h\in P_\gamma$.
Therefore a finite pencil $\Pi^\fin_\gamma$
  sweep all space $\E_n$.

Second, it is easy to check,
that a geodesic $h\gamma(t)h^\top$,
where $h\in P_\gamma$, has no intersections with $\gamma$
or coincide with $\gamma$.
\hfill $\square$

Fix a point $X_0\in\E_n$ (to be concrete, let $X_0=E$).
Denote by $T$ the tangent space to $\E_n$ at $X_0$.
Evidently, we can consider $T$ as the space
of all symmetric matrices. By $\P T$
we denote the set of all  rays in $T$;
a {\it ray} is a set of the form $\theta v$,
where a nonzero vector $v\in T$
is fixed and $\theta$ ranges in positive numbers.

For each ray $\xi$, we consider the geodesic
$\gamma_\xi$ passing through $X_0$ in the direction $\xi$.

By Proposition 4.1, the set of all finite pencils
 is in one-to-one correspondence
with the space $\P T$.

Now we are ready to glue the sphere $S^\far$  at infinity to $\E_n$.
Points $A_\xi$ of the sphere $S^\far$ are enumerated by
 rays $\xi\in \P T$. It remains to define the convergence.

Consider a sequence $Z_1,Z_2,\dots\in \E_n$.
Consider the geodesics $\gamma_{\xi_j}$ connecting $X_0=E$
with $Z_j$. The sequence $Z_j$ converges to a point
$A_\xi$ iff
 $\lim\limits_{j\to\infty}\rho(E,Z_j)=\infty$
and
 $\lim\limits_{j\to\infty} \xi_j=\xi$.

\smallskip

{\bf 4.2. Tilling of $S^\far$.}
Thus we identified the space of finite pencils with
the sphere $S^\far$. Another parametrization
of the same space was given above in Subsection 3.8.
This parametrization gives a canonical tilling
of the sphere $S^\far$
by
a continual family of open
simplexes.

Consider an arbitrary flag  $\W\in \F(I)$, see 2.1.
Consider the set $\Delta(\W)$
of all finite pencils, whose
limit flag
(see Subsection 3.8) is $\W$.
 By 3.8,
$\Delta(\W)\simeq\Delta(I)$.

{\sc Proposition 4.2.}
{\it
The closure of
$\Delta(\W)$
is
$
\bigcup_{\W'\subset \W} \Delta(\W')$
where $\{W'\}$ ranges in all subflags of $\W$.}

{\sc Remark.} This structure is called
 {\it Tits building at infinity}.
Its abstract definition for an  arbitrary
space of nonpositive curvature is contained
in \cite{BGS}.

\medskip

{\bf 5. Solvable pencils.
Karpelevich  and associahedral boundaries}

\medskip

  \addtocounter{sec}{1}
 \setcounter{equation}{0}

In this section we consider the
symmetric spaces $\PE_n$.

{\bf 5.1. The inductive definition of the
 associahedral
 boundary.}
We intend to construct
the associahedral compactification $\Ass(\PE_n)$
of the spaces $\PE_n$ (see \cite{Ner4}).
First, we will describe these compactifications as  disjoint
unions of sets (as it was done above for $\S_n$).

The existence of a natural topology on $\Ass(\PE_n)$
is claimed in Theorem  5.1, the explicit
construction is contained  in Section 7;
before this, in Section 6 we  describe the closure
of a Cartan subspace in the associahedral
compactification.

\smallskip

The construction of the compactification is inductive.
Assume that $\Ass(\PE_k)$ is constructed for all $k<n$.

For any solvable pencil $\Pi_\gamma^\sol$,
we define its limit point at infinity
as the corresponding collection
of the boundary data from Subsection 3.7.

By Corollary 3.7, the boundary obtained in this way is the
union of $2^{n-1}-1$ disjoint pieces
$\K(I)$ having the same dimension
$n^2-2=\dim\PE_n-1.$
Each piece has the form
\begin{equation}
\K(I)=\Delta_\circ(I)\times
\S_n(I)
.\end{equation}

The space $\S_n(I)$
is a bundle, whose base is the space
of (noncomplete) flags $\F_n(I)$ (defined
above in 2.1) and  fibers are the
symmetric spaces
\begin{equation}
\prod\nolimits_s \PE_{i_s-i_{s-1}}
.\end{equation}

We will call these fibers  by
{\it boundary symmetric spaces}.

Now we assume that each fibre (5.2)
already is compactified as
\begin{equation}
\prod\nolimits_s \Ass(\PE_{i_s-i_{s-1}})
,\end{equation}
all the spaces  $\Ass(\PE_{i_s-i_{s-1}})$
are defined by the inductive
hypothesis.

Thus, the boundary is constructed.

{\sc Remark.} We compactified the factor
$\S_n(I)$ in (5.1).
 We emphasis that  the simplexes $\Delta_\circ(I)$
are not compact and hence a topology of a compact
space is yet not defined.

{\bf 5.2. Existence of topology.}

\nopagebreak

{\sc Theorem 5.1.}
{\it There exists a topology of a compact metrizable
space on each $\Ass(\PE_n)$ such that for any
 geodesic $\gamma$ in the space $\PE_n$ or in
any boundary symmetric space the limit
$\lim\nolimits_{t\to+\infty} \gamma(t)$
 with respect to
this topology  coincides
with the limit in the sense defined above.}

{\bf 5.3. Inductive construction
of the Karpelevich boundary.}
Now we intend to construct the Karpelevich
compactification
$\Karp(\PE_n)$ of $\PE_n$.
Its
inductive construction given below
involves   the Karpelevich compactifications
of all the spaces $\prod\PE_{n_j}$.

\smallskip

Thus, assume that the compactifications
$\Karp(\PE_{k_1}\times\dots\times\PE_{k_l})$
are already constructed for
all the collections $(k_1,\dots,\dots, k_l)$
such that
$k_1+\dots+k_l<n.$
Consider a space
$
{\cal E}:=
  \PE_{m_1}\times\dots\times\PE_{m_\beta}$,
where $m_1+\dots+m_\beta=n$.

For each solvable pencil $\Pi_\gamma^\sol$ in $\cal E$,
we add formally a corresponding point at infinity.
Thus
the
set of all such points is a disjoint union of strata
$\K(I_1,\dots,I_\beta)$ described above in
Subsection 3.12.

Each stratum is a product of some polyhedron
$\Delta_\circ(I_1,I_2,\dots)$
and the set
$\S_{m_1}(I_1)\times \dots \times\S_{m_\beta}(I_\beta)$.
The latter set is a bundle,
whose base is
$\prod \F_{m_\tau}(I_\tau)$ and  fibers have
 the form
$$
\prod\limits_{\tau=1}^\beta
\prod\limits_{k=1}^{j_\tau}
\PE_{i_k^{(\tau)}-i_{k-1}^{(\tau)}}
.$$
After this, we replace each fibre by its
compactification
$$
\Karp\Bigl(\prod\limits_{\tau=1}^\beta
\prod\limits_{k=1}^{j_\tau}
\PE_{i_k^{(\tau)}-i_{k-1}^{(\tau)}}\Bigr)
.$$
These space are constructed  by the inductive
hypothesis.

{\bf 5.4. Existence theorem.}

{\sc Theorem 5.2.}
{\it The exists a topology of a
compact metrizable
 space on each set
$
\Karp\bigl( \PE_{m_1}\times\dots\times\PE_{m_\beta}\bigr)
$
and limits of geodesics with respect to this topology
coincide with limits constructed above.}

\medskip

{\bf 6. Permutoassociahedron and Karpelevich polyhedron.}

\medskip

\nopagebreak

  \addtocounter{sec}{1}
 \setcounter{equation}{0}

\nopagebreak

\def\3{
{\sf Fig. 3.}
\begin{picture}(300,135)(0,-45)
\put(40,40){\line(0,1){40}}
\put(40,40){\line(0,-1){40}}
\qbezier(6,20)(40,40)(74,60)
\qbezier(6,60)(40,40)(74,20)
\put(82.5,15){\line(0,1){50}}
\put(85,40){$\Delta_3$}
\put(0,-20){\!\!\!\!\!\!\!\!\!\!\sf \small a) The space $\Xi_3\simeq\R^2$,
the lines $\phi_i=\phi_j$,}
\put(0,-30)
{\!\!\!\!\!\!\!\!\!\!\sf \small and the simplex
$\partial\Lambda$ at infinity}

\put(200,10){
\put(0,0){\line(0,1){40}}
\put(0,0){\line(1,0){40}}
\put(40,0){\line(0,1){40}}
\put(0,40){\line(1,0){40}}
\put(55,55){\line(0,-1){40}}
\put(55,55){\line(-1,0){40}}
\put(40,0){\line(1,1){15}}
\put(40,40){\line(1,1){15}}
\put(0,40){\line(1,1){15}}
\multiput(0,0)(2,2){20}{\circle*{0.2}}
\multiput(40,0)(-2,2){20}{\circle*{0.2}}
\multiput(0,40)(2.75,0.75){20}{\circle*{0.2}}
\multiput(40,0)(0.75,2.75){20}{\circle*{0.2}}
\multiput(15,55)(1.25,-0.75){20}{\circle*{0.2}}
\multiput(55,15)(-0.75,1.25){20}{\circle*{0.2}}
}
\put(170,-20)
{\sf \small b) 24 Weyl chambers in $\Xi_4\simeq\R^3$.}
\put(170,-30)
{\sf \small
Intersections of planes $\phi_i=\phi_j$}
\put(170,-40)
{\sf \small
 and the surface of cube.}
\end{picture}
}

\def\4{\!\!\!\!\!\
\begin{picture}(300,130)(50,0)
\put(-50,0){  
\put(130,0){\small $(1234567)$}
\put(150,10){\line(-2,1){40}}
\put(150,10){\line(3,1){60}}
\put(220,20){$(4567)$}
\put(80,35){\small $(123)$}
\put(110,30){\line (-1,3){20}}
\put(110,30){\line(0,1){60}}
\put(210,30){\line(-1,2){20}}
\put(190,70){\line(-1,1){20}}
\put(190,70){\line(0,1){20}}
\put(170,60){\small $(45)$}
\put(210,30){\line(1,1){20}}
\put(230,50){\line(0,1){40}}
\put(230,50){\line(1,2){20}}
\put(240,55){$\small (67)$}
\put(245,95){$\small (7)$}
\put(225,95){$\small (6)$}
\put(185,95){$\small (5)$}
\put(165,95){$\small (4)$}
\put(105,95){$\small (23)$}
\put(85,95){$\small (1)$}
\multiput(70,10)(5,0){45}{\line(1,0){2}}
\multiput(70,30)(5,0){45}{\line(1,0){2}}
\multiput(70,50)(5,0){45}{\line(1,0){2}}
\multiput(70,70)(5,0){45}{\line(1,0){2}}
\multiput(70,90)(5,0){45}{\line(1,0){2}}
\put(305,95){ \sf Fig. 4}
\put(305,85){\small\sf The tree-partition is (1234567), (123), }
\put(305,75){\small\sf (4567), (23) (45), (67), (1), (4), (5), (6), (7)}
\put(305,65){\small\sf The leveled tree-partition is}
\put(305,55){\small\sf  (1234567);}
\put(305,45){\small\sf (123) (4567);}
\put(305,35){\small\sf (1) (23) (4567);}
\put(305,25){\small\sf (1) (23) (45) (67);}
\put(305,15){\small\sf (1) (23) (45) (6) (7);}
\put(305,5){\small\sf (1) (23) (4) (5) (6) (7).}
}
\end{picture}
}

\def\5{
\begin{picture}(250,50)
\put(0,10){
\put(80,0){\line(-2,1){80}}
\put(80,0){\line(2,1){80}}
\put(80,0){\line(1,1){40}}
\put(80,0){\line(-1,1){40}}
\put(80,0){\line(0,1){40}}
\put(80,0){\line(0,-1){10}}

\put(0,30){$a_1$}
\put(33,30){$a_2$}
\put(70,30){$a_3$}
\put(100,30){$a_4$}
\put(130,30){$a_5$}

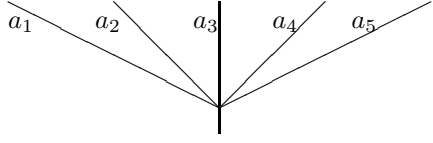
\put(170,35){\small\sf Fig.5 A churn-staff of $\frA$.The space $U_\frA(K)$ }
\put(170,25){\small\sf consists of collections $(a_1:a_2:a_3:a_4:a_5)$}
\put(170,15){\small\sf defined up to positive factor and addition of}
\put(170,5){\small\sf  constant. Numbers $a_j$ are pairwise different}
}
\end{picture}
}

\def\6{{\smallskip
\begin{picture}(250,50)
\multiput(10,10)(5,0){45}{\line(1,0){2}}
\multiput(10,50)(5,0){45}{\line(1,0){2}}

\put(30,10){\line(-1,1){40}}
\put(30,10){\line(0,1){40}}
\put(30,10){\line(1,1){40}}
\put(30,10){\line(0,-1){10}}

\put(-5,30){$a_1$}
\put(18,30){$a_2$}
\put(40,30){$a_3$}

\put(90,10){\line(0,1){40}}
\put(90,10){\line(0,-1){10}}
\put(80,30){$b_1$}

\put(200,10){\line(-1,1){40}}
\put(200,10){\line(0,1){40}}
\put(200,10){\line(1,1){40}}
\put(200,10){\line(0,-1){10}}

\put(165,30){$c_1$}
\put(188,30){$c_2$}
\put(210,30){$c_3$}

\end{picture}
\sf\small

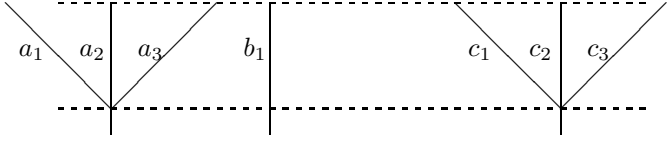
Fig 6. The set $W_j[\frA]$ consists of collections
$u=(a_1:a_2:a_3:b_1:c_1:c_2:c_3:\dots)$
defined up to a  common {\it\small positive} factor and up to an addition
of $(t:t:t\,:\,s\,:\,r:r:r\,:\dots)$.
The numbers $a_1$, $a_2$, $a_3$
are pairwise different;
$b_1$, $b_2$, $b_3$  are pairwise different, etc.
 In our case the variable $b_1$
is fake.

  The map $\Karp_n\to\Pass_n$ takes $u$ to
the collection $[(a_1:a_2:a_3), (c_1:c_2:c_3),\dots]$.
In each bracket (\dots), the numbers are defined
up to a common positive factor and addition of $(t:t:\dots)$. }
}

\def\7{
{
\setlength{\unitlength}{2pt}
\begin{picture}(130,160)(10,-50)
\put(-10,-10){\line(1,1){110}}
\put(-10,0){\line(1,0){110}}
\put(90,-10){\line(0,1){110}}
\put(80,90){\line(1,0){20}}
\put(0,10){\line(0,-1){20}}
{\thicklines
\put(70,10){\line(1,1){10}}
\put(80,20){\line(0,1){20}}
\put(80,40){\line(-2,1){20}}
\put(60,50){\line(-1,-1){20}}
\put(40,30){\line(1,-2){10}}
\put(50,10){\line(1,0){20}}
\put(70,10){\line(0,-1){20}}
\put(80,20){\line(1,0){20}}
\put(80,40){\line(1,0){20}}
\put(60,50){\line(-1,1){10}}
\put(40,30){\line(-1,1){10}}
\put(50,10){\line(0,-1){20}}
}
\multiput(70,-10)(1,-1){10}{\circle*{0.1}}
\multiput(80,-20)(1,0){20}{\circle*{0.1}}
\multiput(100,-20)(1,1){10}{\circle*{0.1}}
\multiput(110,-10)(0,1){20}{\circle*{0.1}}
\multiput(110,10)(-1,1){10}{\circle*{0.1}}
\multiput(100,20)(0,1){20}{\circle*{0.1}}
\multiput(100,40)(1,0.5){20}{\circle*{0.1}}
\multiput(50,60)(-0.5,1){20}{\circle*{0.1}}
\multiput(50,60)(-1,-1){20}{\circle*{0.1}}
\multiput(30,40)(-1,0.5){20}{\circle*{0.1}}
\multiput(50,-10)(-0.5,-1){20}{\circle*{0.1}}
\multiput(50,-10)(1,0){20}{\circle*{0.1}}
\put(70,12){\line(1,1){8}}
\put(72,8){\small $A$}
\put(72,-10){\small $A'$}
\put(80,15){\small $B$}
\put(95,15){\small $B'$}
\put(80,42){\small $C$}
\put(95,42){\small $C'$}
\put(61,51){\small $D$}
\put(51,61){\small $D'$}
\put(35,25){\small $E$}
\put(25,35){\small $E'$}
\put(45,8){\small $F$}
\put(45,-10){\small $F'$}
\put(2,-5){\small $P$}
\put(92,-5){\small $Q$}
\put(92,85){\small $T$}
\put(0,-30){\small Fig.7}
\end{picture}}
}

\def\8{{
{\sc Fig.8.} \small The list of faces in a given Weyl chamber.

1) The unique 3-dimensional cell. On Fig. 7
 it is under the sheet of the paper.
{\setlength{\unitlength}{2pt}
\begin{picture}(40,10)(0,0)
\put(10,0){\circle*{2}}
\put(3,5){\small $(1234)$}
\end{picture}}

2) 2-dimensional faces

{\setlength{\unitlength}{2pt}
\begin{picture}(150,35)(0,0)
\put(20,10){\line(-1,1){10}}
\put(20,10){\line(1,1){10}}
\put(10,20){\circle*{2}}
\put(30,20){\circle*{2}}
\put(20,10){\circle*{2}}
\put(5,25){\small $(123)$}
\put(25,25){\small $(4)$}
\put(5,0){\small $\dots E'EFF'\dots$}
\put(35,0){
\put(20,10){\line(-1,1){10}}
\put(20,10){\line(1,1){10}}
\put(10,20){\circle*{2}}
\put(30,20){\circle*{2}}
\put(20,10){\circle*{2}}
\put(5,25){\small $(12)$}
\put(25,25){\small $(34)$}
\put(5,0){\small $\dots A'ABB'\dots$}
}
\put(70,0){
\put(20,10){\line(-1,1){10}}
\put(20,10){\line(1,1){10}}
\put(10,20){\circle*{2}}
\put(30,20){\circle*{2}}
\put(20,10){\circle*{2}}
\put(5,25){\small $(1)$}
\put(25,25){\small $(234)$}
\put(5,0){\small $\dots C'CDD'\dots$}
}
\put(120,0){
\put(20,10){\line(-1,1){10}}
\put(20,10){\line(1,1){10}}
\put(20,10){\line(2,1){20}}
\put(20,10){\line(-2,1){20}}
\put(10,20){\circle*{2}}
\put(30,20){\circle*{2}}
\put(20,10){\circle*{2}}
\put(0,20){\circle*{2}}
\put(40,20){\circle*{2}}
\put(-3,25){\small $(1)$}
\put(7,25){\small $(2)$}
\put(27,25){\small $(3)$}
\put(37,25){\small $(4)$}
\put(5,0){\small $ABCDEF$}
}
\end{picture}}

{\setlength{\unitlength}{2pt}
\begin{picture}
(150,35)(0,0)
\put(20,10){\line(-1,1){10}}
\put(20,10){\line(0,1){10}}
\put(20,10){\line(1,1){10}}
\put(10,20){\circle*{2}}
\put(30,20){\circle*{2}}
\put(20,10){\circle*{2}}
\put(20,20){\circle*{2}}
\put(5,25){\small $(12)$}
\put(17,25){\small $(3)$}
\put(27,25){\small $(4)$}
\put(5,0){\small $ F'FAA'$}
\put(35,0)
{
\put(20,10){\line(-1,1){10}}
\put(20,10){\line(0,1){10}}
\put(20,10){\line(1,1){10}}
\put(10,20){\circle*{2}}
\put(30,20){\circle*{2}}
\put(20,10){\circle*{2}}
\put(20,20){\circle*{2}}
\put(7,25){\small $(1)$}
\put(16,25){\small $(23)$}
\put(27,25){\small $(4)$}
\put(5,0){\small $D'DEE'$}
}
\put(70,0)
{
\put(20,10){\line(-1,1){10}}
\put(20,10){\line(0,1){10}}
\put(20,10){\line(1,1){10}}
\put(10,20){\circle*{2}}
\put(30,20){\circle*{2}}
\put(20,10){\circle*{2}}
\put(20,20){\circle*{2}}
\put(7,25){\small $(1)$}
\put(16,25){\small $(2)$}
\put(27,25){\small $(34)$}
\put(5,0){\small $ B'BCC'$}
}

\end{picture}}

3) Edges

\begin{picture}(150,50)

\put(10,20){\line(1,1){10}}
\put(20,10){\line(-1,1){20}}
\put(20,10){\line(1,1){20}}
\put(20,10){\circle*{2}}
\put(0,30){\circle*{2}}
\put(40,30){\circle*{2}}
\put(20,30){\circle*{2}}
\put(10,20){\circle*{2}}
\put(-5,35){\small $(12)$}
\put(17,35){\small $(3)$}
\put(37,35){\small $(4)$}
\put(15,0){\small $FF'$}

\put(65,0)
{
\put(10,20){\line(1,1){10}}
\put(20,10){\line(-1,1){20}}
\put(20,10){\line(1,1){20}}
\put(20,10){\line(1,1){20}}
\put(20,10){\circle*{2}}
\put(20,10){\circle*{2}}
\put(0,30){\circle*{2}}
\put(40,30){\circle*{2}}
\put(10,20){\circle*{2}}
\put(20,30){\circle*{2}}
\put(-5,35){\small $(1)$}
\put(15,35){\small $(23)$}
\put(37,35){\small $(4)$}
\put(15,0){\small $EE'$}
}

\put(130,0)
{
\put(30,20){\line(-1,1){10}}
\put(20,10){\line(-1,1){20}}
\put(20,10){\line(1,1){20}}
\put(30,20){\circle*{2}}
\put(0,30){\circle*{2}}
\put(40,30){\circle*{2}}
\put(20,30){\circle*{2}}
\put(-5,35){\small $(1)$}
\put(15,35){\small $(23)$}
\put(37,35){\small $(4)$}
\put(15,0){\small $DD'$}
}

\put(195,0)
{
\put(30,20){\line(-1,1){10}}
\put(20,10){\line(-1,1){20}}
\put(20,10){\line(1,1){20}}
\put(30,20){\circle*{2}}
\put(0,30){\circle*{2}}
\put(40,30){\circle*{2}}
\put(20,30){\circle*{2}}
\put(-5,35){\small $(1)$}
\put(15,35){\small $(2)$}
\put(37,35){\small $(34)$}
\put(15,0){\small $CC'$}
}

\put(260,0)
{
\put(30,20){\line(-1,1){10}}
\put(20,10){\line(-1,1){20}}
\put(20,10){\line(1,1){20}}
\put(30,20){\circle*{2}}
\put(0,30){\circle*{2}}
\put(40,30){\circle*{2}}
\put(20,30){\circle*{2}}
\put(-5,35){\small $(12)$}
\put(15,35){\small $(3)$}
\put(37,35){\small $(4)$}
\put(15,0){\small $AA'$}
}

\end{picture}

\begin{picture}(150,50)

\put(10,20){\line(1,1){10}}
\put(20,10){\line(-1,1){20}}
\put(20,10){\line(1,1){20}}
\put(20,10){\circle*{2}}
\put(0,30){\circle*{2}}
\put(40,30){\circle*{2}}
\put(20,30){\circle*{2}}
\put(10,20){\circle*{2}}
\put(-5,35){\small $(1)$}
\put(17,35){\small $(2)$}
\put(37,35){\small $(34)$}
\put(15,0){\small $BB'$}

\put(70,0){
\put(20,10){\line(-1,1){20}}
\put(20,10){\line(1,1){20}}
\put(10,20){\line(0,1){10}}
\put(30,20){\line(0,1){10}}
\multiput(0,20)(2,0){20}{\circle*{0.5}}
\put(15,0){\small $AB$}
\put(-5,35){\small $(1)$}
\put(6,35){\small $(2)$}
\put(35,35){\small $(4)$}
\put(24,35){\small $(3)$}
\put(20,10){\circle*{2}}
\put(10,20){\circle*{2}}
\put(30,20){\circle*{2}}
\put(10,30){\circle*{2}}
\put(30,30){\circle*{2}}
\put(0,30){\circle*{2}}
\put(40,30){\circle*{2}}
}

\put(140,0){
\put(20,10){\line(-1,1){20}}
\put(20,10){\line(1,1){20}}
\put(20,10){\line(0,1){20}}
\put(10,20){\line(0,1){10}}
\put(15,0){\small $AF$}
\put(-5,35){\small $(1)$}
\put(6,35){\small $(2)$}
\put(35,35){\small $(4)$}
\put(18,35){\small $(3)$}
\put(20,10){\circle*{2}}
\put(10,20){\circle*{2}}
\put(10,30){\circle*{2}}
\put(20,30){\circle*{2}}
\put(0,30){\circle*{2}}
\put(40,30){\circle*{2}}
}

\put(210,0)
{
\put(20,10){\line(-1,1){20}}
\put(20,10){\line(1,1){20}}
\put(20,10){\line(0,1){10}}
\put(20,20){\line(1,1){10}}
\put(20,20){\line(-1,1){10}}
\put(20,10){\circle*{2}}
\put(10,30){\circle*{2}}
\put(20,20){\circle*{2}}
\put(30,30){\circle*{2}}
\put(0,30){\circle*{2}}
\put(40,30){\circle*{2}}
\put(14,0){\small $DE$}
\put(-5,35){\small $(1)$}
\put(6,35){\small $(2)$}
\put(35,35){\small $(4)$}
\put(23,35){\small $(3)$}
}

\put(280,0){
\put(20,10){\line(-1,1){20}}
\put(20,10){\line(1,1){20}}
\put(20,10){\line(0,1){20}}
\put(30,20){\line(0,1){10}}
\put(15,0){\small $BC$}
\put(-5,35){\small $(1)$}
\put(16,35){\small $(2)$}
\put(38,35){\small $(4)$}
\put(27,35){\small $(3)$}
\put(20,10){\circle*{2}}
\put(20,30){\circle*{2}}
\put(30,30){\circle*{2}}
\put(0,30){\circle*{2}}
\put(40,30){\circle*{2}}

}

\end{picture}

\begin{picture}(150,50)

\put(10,20){\line(1,1){10}}
\put(10,20){\line(0,1){10}}
\put(20,10){\line(-1,1){20}}
\put(20,10){\line(1,1){20}}
\put(20,10){\circle*{2}}
\put(0,30){\circle*{2}}

\put(40,30){\circle*{2}}
\put(20,30){\circle*{2}}
\put(10,20){\circle*{2}}
\put(10,30){\circle*{2}}
\put(-5,35){\small $(1)$}
\put(7,35){\small $(2)$}
\put(18,35){\small $(3)$}
\put(37,35){\small $(4)$}
\put(15,0){\small $EF$}

\put(80,0){
\put(20,10){\line(-1,1){20}}
\put(20,10){\line(1,1){20}}
\put(30,20){\line(-1,1){10}}
\put(30,20){\line(0,1){10}}
\put(15,0){\small $CD$}
\put(-5,35){\small $(1)$}
\put(16,35){\small $(2)$}
\put(38,35){\small $(4)$}
\put(27,35){\small $(3)$}
\put(20,10){\circle*{2}}
\put(20,30){\circle*{2}}
\put(30,30){\circle*{2}}
\put(0,30){\circle*{2}}
\put(40,30){\circle*{2}}

}

\end{picture}

4) Vertices

\begin{picture}(150,70)(0,0)

\put(80,0){
\put(30,10){\line(-1,1){30}}
\put(30,10){\line(1,1){30}}
\put(20,20){\line(0,1){20}}
\put(50,30){\line(0,1){10}}
\multiput(0,20)(2,0){30}{\circle*{0.5}}
\multiput(0,30)(2,0){30}{\circle*{0.5}}
\put(30,10){\circle*{2}}
\put(20,20){\circle*{2}}
\put(50,30){\circle*{2}}
\put(0,40){\circle*{2}}
\put(60,40){\circle*{2}}
\put(20,40){\circle*{2}}
\put(50,40){\circle*{2}}
\put(-5,45){$(1)$}
\put(15,45){$(2)$}
\put(43,45){$(3)$}
\put(55,45){$(4)$}
\put(25,0){$B$}
}

\put(30,10){\line(-1,1){30}}
\put(30,10){\line(1,1){30}}
\put(40,20){\line(0,1){20}}
\put(10,30){\line(0,1){10}}
\multiput(0,20)(2,0){30}{\circle*{0.5}}
\multiput(0,30)(2,0){30}{\circle*{0.5}}
\put(30,10){\circle*{2}}
\put(40,20){\circle*{2}}
\put(10,30){\circle*{2}}
\put(0,40){\circle*{2}}
\put(60,40){\circle*{2}}
\put(10,40){\circle*{2}}
\put(40,40){\circle*{2}}
\put(-5,45){$(1)$}
\put(6,45){$(2)$}
\put(33,45){$(3)$}
\put(55,45){$(4)$}
\put(25,0){$A$}

\put(160,0)
{
\put(30,10){\line(-1,1){30}}
\put(30,10){\line(1,1){30}}
\put(20,20){\line(1,1){20}}
\put(10,30){\line(1,1){10}}
\multiput(0,20)(2,0){30}{\circle*{0.5}}
\multiput(0,30)(2,0){30}{\circle*{0.5}}
\put(30,10){\circle*{2}}
\put(20,20){\circle*{2}}
\put(10,30){\circle*{2}}
\put(0,40){\circle*{2}}
\put(60,40){\circle*{2}}
\put(20,40){\circle*{2}}
\put(40,40){\circle*{2}}
\put(-5,45){$(1)$}
\put(16,45){$(2)$}
\put(33,45){$(3)$}
\put(55,45){$(4)$}
\put(25,0){$F$}
}

\put(240,0)
{
\put(30,10){\line(-1,1){30}}
\put(30,10){\line(1,1){30}}
\put(20,20){\line(1,1){20}}
\put(30,30){\line(-1,1){10}}
\multiput(0,20)(2,0){30}{\circle*{0.5}}
\multiput(0,30)(2,0){30}{\circle*{0.5}}
\put(30,10){\circle*{2}}
\put(20,20){\circle*{2}}
\put(30,30){\circle*{2}}
\put(0,40){\circle*{2}}
\put(60,40){\circle*{2}}
\put(20,40){\circle*{2}}
\put(40,40){\circle*{2}}
\put(-5,45){$(1)$}
\put(16,45){$(2)$}
\put(33,45){$(3)$}
\put(55,45){$(4)$}
\put(25,0){$E$}
}

\end{picture}

\begin{picture}(150,70)(0,0)

\put(30,10){\line(-1,1){30}}
\put(30,10){\line(1,1){30}}
\put(40,20){\line(-1,1){20}}
\put(30,30){\line(1,1){10}}
\multiput(0,20)(2,0){30}{\circle*{0.5}}
\multiput(0,30)(2,0){30}{\circle*{0.5}}
\put(30,10){\circle*{2}}
\put(40,20){\circle*{2}}
\put(30,30){\circle*{2}}
\put(0,40){\circle*{2}}
\put(60,40){\circle*{2}}
\put(20,40){\circle*{2}}
\put(40,40){\circle*{2}}
\put(-5,45){$(1)$}
\put(16,45){$(2)$}
\put(33,45){$(3)$}
\put(55,45){$(4)$}
\put(25,0){$D$}

\put(90,0){
\put(30,10){\line(-1,1){30}}
\put(30,10){\line(1,1){30}}
\put(40,20){\line(-1,1){20}}
\put(50,30){\line(-1,1){10}}
\multiput(0,20)(2,0){30}{\circle*{0.5}}
\multiput(0,30)(2,0){30}{\circle*{0.5}}
\put(30,10){\circle*{2}}
\put(40,20){\circle*{2}}
\put(50,30){\circle*{2}}
\put(0,40){\circle*{2}}
\put(60,40){\circle*{2}}
\put(20,40){\circle*{2}}
\put(40,40){\circle*{2}}
\put(-5,45){$(1)$}
\put(16,45){$(2)$}
\put(33,45){$(3)$}
\put(55,45){$(4)$}
\put(25,0){$C$}
}

\end{picture}
}}

\nopagebreak

{\bf A. Definition of permutoassociahedron and
karpelevich-hedron}

\smallskip

\nopagebreak


{\bf 6.1. Spaces $\Xi(I)$.}
Let $I$ be a finite set,
denote by $\#(I)$ the number of its elements;
 the basic example is
the set $I=\bJ$:
$$\bJ=\bJ_n:=\{1,2,\dots,n\}.$$
Denote by $\Xi(I)$ the set of all functions
$I\to \R$ defined up to an addition of a constant
function;
we also denote
$$\Xi_n:=\Xi(\bJ).$$
The  space $\Xi_n$
consists of
vectors $(\phi_1,\dots,\phi_{n})\in\R^n$
 defined up
to the equivalence
$$
(\phi_1,\dots,\phi_n)\sim
 (\phi_1+a,\dots,\phi_n+a)
.$$
We also can consider elements of $\Xi_n$
as ordered collections of points on $\R$
defined up to a translation.

\smallskip

For a subset $K\subset I$,
we have the natural map
\begin{equation}
\Xi(I)\to\Xi(K)
\end{equation}
(we forget part of coordinates).

Consider a partition $\fra$ of the set $I$,
denote by $I/\fra$ the corresponding quotient
and by $\pi:I\to I/\fra$ the natural
projection.
We have a natural embedding
\begin{equation}
\Xi(I/\fra)\to\Xi(I)
,\end{equation}
i.e., to a function $f:I/\fra\to\R$
we assign the function
$f\circ\pi$.

Again, consider a partition $\fra$ of $I$, let
$I_1$, \dots, $I_s$ be its elements.
Denote by $C[I;\fra]$ the space of functions
$I\to\R$ that are constants on each subset
$I_m$.
Consider the quotient linear space
\begin{equation}
\Xi[I;\fra]:=\Xi(I)/C[I;\fra]
.\end{equation}

We have natural projection map
\begin{equation}
\Xi(I)\to \Xi[I;\fra]
.\end{equation}
Also, we have the obvious identification
\begin{equation}
\Xi[I;\fra]\simeq \bigoplus\nolimits_{m=1}^s \Xi(I_m)
.\end{equation}

If a partition $\frb$ is a subdivision of $\fra$,
then we have the map
$$
\qquad\qquad \qquad  \qquad \qquad
\Xi[I;\fra]\to \Xi[I;\frb]
\qquad\qquad \qquad  \qquad \qquad \qquad \qquad
(6.5.a)              \!    \!\!\!
$$

We call by {\it walls}
the hyperplanes $f(a)=f(b)$,
where $a,b\in I$. These hyperplanes divide the
space $\Xi(I)$ into $\#(I)!$ simplicial cones,
which are called {\it Weyl chambers}.

{\sc Example.}
  Tilings of $\Xi_3$, $\Xi_4$ by Weyl chambers are presented on Fig. 3.

\smallskip

Now assume that the set $I$ is an ordered
set with the order $\prec$.
Then we have the {\it positive
Weyl chamber} $\Lambda^+(I)$ defined by
the inequalities
$$
a\prec b \quad \Rightarrow  \quad f(a)\ge f(b)
$$
We also denote
$\Lambda^+_n:=\Lambda^+(\bJ)$.

\begin{figure}
\3
\end{figure}

{\bf 6.2. Compactification of the spaces $\Xi(I)$.}
Let $V$ be a linear space.
A {\it ray} is a subset in $V$ having the form
$\lambda v$, where $v\ne0$ is a fixed vector in $V$
and $\lambda$ ranges in positive numbers.

We compactify each ray by a point at infinity.
Denote the set of all such points at infinity
by $\partial V$ ({\it sphere at infinity}).
By $\overline V$ we denote $V\cup\partial V$.
We define a topology on $\overline V$
in the obvious way.

In particular, we obtain the spaces
$$\overline{\Xi(I)}=\Xi(I)\cup\partial\Xi(I),
\qquad
\overline{\Xi[I;\fra]}=
\Xi[I;\fra]\cup \partial\Xi[I;\fra]
.$$
We emphasis (compare with (6.5)), that
$$
\overline{\Xi[I;\fra]}\ne \prod \overline{\Xi(I_k)}
.$$

For an ordered set $I$ we denote
by $\overline{\Lambda^+(I)}$ the closure
of the positive Weyl chamber $\Lambda^+(I)$ in $\overline{\Xi(I)}$
and  by
$\partial \Lambda^+(I):=\overline{\Lambda^+(I)}\setminus\Lambda^+(I)$
its boundary.

{\it This allows to consider $\overline{\Xi(I)}$
 as a polyhedron}; the space $\Xi(I)$
is its interior and
the boundary $\partial \Xi(I)$
is divided into simplexes of the type $\partial\Lambda^+$.
This point of view is also represented on Fig. 3.

{\bf 6.3. Definition of the permutoassociahedron.}
For each subset $I\subset \bJ$
 consider the 'forgetting' map
$\Xi_n\to \Xi(I)$, see (6.1).
Consider the diagonal embedding
$$
\iota_n:\Xi_n\to\prod\limits_{I\subset\bJ} \Xi(I)
.$$
We also have the inclusion
\begin{equation}
\prod\limits_{I\subset\bJ} \Xi(I)
\subset
\prod\limits_{I\subset\bJ} \overline\Xi(I)
.\end{equation}

{\sc Definition.}
The permutoassociahedron $\Pass_n$ (see \cite{Kap})
 is the closure of
image
$\iota_n(\Xi_n)$ in
$\prod\nolimits_{I\subset\bJ}
 \overline\Xi(I)$.

{\bf 6.4. Definition of the Karpelevich polyhedron.}
For each partition $\fra$ of $\bJ$,
we have the  map
$\Xi_n\to \Xi[\bJ;\fra]$, see (6.4).
Consider the diagonal embedding
\begin{equation}
\Xi_n\to\prod\limits_{\fra} \Xi[\bJ;\fra]\subset
\prod\limits_{\fra} \overline{\Xi[\bJ;\fra]}
,\end{equation}
where the product is given over all the partitions $\fra$ of
$\bJ$.

{\sc Definition.}
The Karpelevich polyhedron $\Karp_n$ is the closure
of the image of $\Xi_n$ in the space
$\prod \overline{\Xi[\bJ;\fra]}$.

{\sc Remark.}
Assume that $\fra$ consists of a subset $I$ and
single-element sets. Then
$\Xi[\bJ;\fra]=\Xi(I)$.
Thus each factor of the product (6.6)
is a factor of the product (6.7), and hence we obtain
the natural projection $\Karp_n\to\Pass_n$.

\smallskip

{\bf B. More notation}

\smallskip

\nopagebreak

{\bf 6.5. Some functorial properties of
spheres at infinity.} This  is used only in
Subsections 6.9 and 6.15.

\nopagebreak

1) Fix a subset $K\subset I$. Denote by
$C(K)$ the space of functions on $I$
which are constants on $K$.
 The map (6.1)
induces the continuous  map
\begin{equation}
\partial {\Xi(I)}\setminus \partial {C(K)}
\to\partial{\Xi(K)}
.\end{equation}

2)
Let $\fra$ be a partition of $I$ with elements $I_k$.
The map (6.2) induces the embedding
\begin{equation}
\overline{\Xi(I/\fra)}\to \overline{\Xi(I)}
.\end{equation}

3) Let $C[I;\fra]$ be the same as in 6.1.
The map (6.4) induces the continuous map
$$
\partial{\Xi(I)}\setminus \partial C[I;\fra]
\to\overline{\Xi[I;\fra]}
.$$

4) Let $\fra$ be a partition of $I$, let $\frb$
be a subdivision of $\fra$. The
map (6.5.a)
induces   the map
$$
\qquad\qquad\qquad\qquad
\partial\Xi[I;\fra] \setminus \partial  C[I;\frb]
\to
\partial\Xi[I;\frb]
\qquad\qquad\qquad\qquad\qquad (6.9.a)  \!\!\!\!\!\!
$$

{\bf 6.6. Notation for sphere at infinity outside walls.}
For each wall $f(a)=f(b)$ in $\Xi(I)$ denote
by $S_{a,b}$ its intersection with $\partial\Xi(I)$.
The $(\#(I)-3)$-dimensional  spheres $S_{a,b}$ divide the
 $(\#(I)-2)$-dimensional sphere $\partial \Xi(I)$ into
 $\#(I)!$ simplexes.
Denote
$$
\partial\Xi(I)^\gen:=
\partial \Xi(I)\setminus \cup S_{a,b}
;$$
see Fig. 3b, it is the surface of cube 
without edges and diagonals of faces.

Also for the Weyl chamber $\Lambda^+(I)$
we denote by
$$\partial\Lambda^+(I)^\gen:=
  \partial\Lambda^+(I)\cap\partial\Xi(I)^\gen
.$$

Now, let $\fra$ be a partition of $I$, let $I_k$ be its elements.
We say that a ray $tf$, where $t>0$, $f\in  \Xi[I;\fra]$
is generic, if for each $I_k$ and
 each $a,b\in I_k$ we have
$f(a)\ne f(b)$.
We  define
 the  set
$\partial \Xi[I;\fra]^\gen\subset\partial \Xi[I;\fra]$
as the set of limits of generic rays.

{\bf 6.7. Combinatorial partition-like
 structures.}

\nopagebreak

{\it Partitions.}
Consider a finite set $\bM  $. Its {\it partition} is a representation
of $\bM  $ as a disjoint union of subsets.


{\it Tree-partitions.}
A system $\frA$ of subsets of $\bM $ is a {\it tree-partition}
if the following conditions are hold

a) $\bM \in\frA$

b) For $I_1,I_2\in \frA$, we have either
$I_1\cap I_2=\emptyset$, or $I_1\supset I_2$,
$I_1\subset I_2$.

c) Let $I\supset K$ be elements of $\frA$. Then there exists
a collection $K_1=K$, $K_2$, \dots, $K_\alpha\in\frA$ such that
\begin{equation}
I=\cup K_j, \qquad \qquad \text{$K_i\cap K_j=\emptyset$
 for $i\ne j$}
\end{equation}

\smallskip

A subset $I\in \frA$
is {\it irreducible}, if there
 is no  $K\in \frA$  such that $K\subset I$.

For a reducible subset $I\in\frA $ there exists its
unique {\it minimal decomposition}
(6.10) such that for  $L\in \frA$
satisfying
$I\supset L\supset K_j$
we have $L=I$ or $L=K_j$.

{\it Another definition of  tree-partitions.}
Consider a set $\bM $. Consider a partition $\frx$ of $\bM $.
For each element  $K_j\in\frx$, consider a partition
$\fry_j$ of $K_j$. Then we repeat the same
with elements of partitions $\fry_j$, etc.
Obviously, we obtain a tree-partition of $\bM$.

{\it Leveled tree-partitions.}
Consider a finite set $\bM$. Its leveled
tree-partition $\frA$
is a family of partitions
\begin{equation}
\fra_0,\fra_1,\dots,\fra_k
\end{equation}
satisfying the conditions

a) $\fra_0$ consists of the set $\bM$ itself.

b)  $\fra_{m+1}$ is a subdivision
of $\fra_m$

c) $\fra_{m+1}\ne \fra_m$ for all $m$.

\smallskip

For a leveled tree-partition $\frA$ of $\bM$ consider
$\cup_m \fra_m$ (i.e., we consider all the elements
of all the partitions $\fra_m$). Obviously,
we obtain a tree-partition of $\bM$.

\begin{figure}
\4
\end{figure}

{\it Segmental partitions.}
Let $\bM$ be an ordered set.
{\it Segments}  $[a,b]\subset \bM$ are subsets
having the form
$a\prec j\prec b$.
A segmental partition (tree-partition, leveled tree-partition)
 is a partition, all whose elements are
 segments.

{\it Perfect tree-partitions.}
A tree-partition is perfect if all its
irreducible elements
are  singletons.

\smallskip

{\bf C. Description of permutoassociahedron}

\smallskip

\nopagebreak

{\bf 6.8. Stratification of the permutoassociahedron.}
The permutoassociahedron $\Pass_n$ was defined as a subset
in the     polyhedron $\prod \overline{\Xi(I)}$,
see (6.6). Considering the intersections of $\Pass_n$
with faces of $\prod \overline{\Xi(I)}$, we obtain
a natural stratification of $\Pass_n$.

Fix a tree-partition $\frA$ of $\bJ$.
First, for any element $K\in\frA$, we intend to define
a set $U_\frA(K)$

a) For an  irreducible $K$, we assume
$U_\frA(K):=\Xi(K)$.

b) Let $K$ be reducible. Let
$\frr$
be its minimal decomposition, and $h(K)$ be number
of its elements.
Then
$U_\frA(K):=\partial \Xi(K/\frr)^\gen$, see Fig.5.

\begin{figure}
\5
\end{figure}

{\sc Remark.} In the reducible case,
 the set $U_\frA(K)$ is a union of $h(K)!$ of disjoint
 $(h(K)-2)$-dimensional
open simplexes. If $h(K)=2$, then
$U_\frA(K)$   is a two-point set.

Now we define the stratum $\Flat(\frA)$
as the product
\begin{equation}
\Flat(\frA)=\Flat^\Pass(\frA):=\prod\nolimits_{K\in\frA} U_\frA(K)
.\end{equation}

{\sc Theorem 6.1.}
{\it The permutoassociahedron is
\begin{equation}
\Pass_n=\bigcup_{\frA} \Flat(\frA)
,\end{equation}
where the union is given over
all the tree-partitions $\frA$ of $\bJ$.}

We emphasis, that a set $\Flat(\frA)$ is disconnected;
 the  number of its components is
$\prod
 h(K)!
,$
the product is given other all reducible elements of $\frA$.
These components are (open) faces of the polyhedron
$\Pass_n$.

{\bf 6.9. Identification of definitions 6.3 and 6.8 of
permutoassociahedron.}
It is sufficient to write a map
\begin{equation}
\Flat(\frA)\to \overline{\Xi(L)}
\end{equation}
for a given tree-partition $\frA$ of the set $\bJ$
and for any subset $L\subset\bJ$.
This will define the canonical
map from (6.13) to $\Pass_n$.

Denote by $K$ the minimal element of
the tree-partition $\fra$ containing $L$.
Obviously, this element exists,
since $\bJ\in\frA$. The image of a point
${ u}\in\Flat(\frA)$
under (6.14) will be completely determined by
its projection to the factor $U_\frA(K)$ in (6.12).
There are two cases: $K\in\frA$ is  irreducible
and
$K\in\frA$ is  reducible.

 First, let $K$ be an irreducible element
of $\frA$. Then
$U_\frA(K)=\overline{\Xi(K)}\to \overline{\Xi(L)}$
is the canonical map (6.1).

Second, let $K$ be reducible. Let $\frr$
be its canonical decomposition.
Then our map is the composition
of the canonical maps (see (6.9), (6.8))
$$
\partial\Xi(K/\frr)^\gen
\to \partial\Xi(K)\to\partial \Xi(L)
.$$

The required map is  constructed.

{\bf 6.10. Convergence
in the permutoassociahedron.}
 Consider a sequence
$x_1$, $x_2, \dots$
in $\Xi_n$.

The first  necessary condition
of the convergence
is the convergence in $\overline\Xi_n$.
If the limit belings to $\Xi_n$, then
it is the limit in $\Pass$
(the corresponding tree-partition consists of
one set $\bJ$).

Otherwise,
let $t(\mu_1,\dots,\mu_n)$ be the limit ray.
We construct a partition $\frp$ of $\bJ$
by the following equivalence relation
\begin{equation}
 k\sim l\qquad \text{if and only if $\mu_k=\mu_l$}
.\end{equation}
Denote the elements of the partition $\frp$ by
$I_1$, $I_2$, \dots. The sequence $x_j$ induces
sequences $x_j^s$ in each space $\Xi(I_s)$.

Our next necessary condition
is the convergence of each sequence
$x_j^s$ in each $\overline{\Xi(I_s)}$,
etc. etc.

{\sc Example.}
Consider the sequence in $\Xi_6$ given by
\begin{equation}
(n^3+2n,
\,\,n^3+n,
\,\,n^3,
\,\,3n,
\,\,2n+1,
\,\,2n)
\end{equation}

Its limit is contained in the set
$\Flat(\frA)$ for the tree-partition
$$\frA:\bigl((1)(2)(3)\bigr)\,\,\,\,  \bigl( (4)\,\, (56)\bigr)$$
Indeed, the limit ray for (6.16) in $\Xi_6$ is
$t(1,1,1,0,0,0)$. This gives the partition
$
(123)(456)
$.

In $\Xi\{1,2,3\}$, we have the sequence
$
(n^3+2n,
\,\,n^3+n,
\,\,n^3)\sim
(2n,n,0)
$.
 Its limit ray is $t(1,1/2,0)\in\partial\Xi\{1,2,3\}$.
This gives the partition $(1)(2)(3)$ of $(123)$.

In $\Xi\{4,5,6\}$, we have the sequence
$(3n,2n+1,2n)\sim (n,1,0)$.
Its limit ray is
$t(1,0,0)\in\partial\Xi\{4,5,6\}$,
and this gives the subpartition of $(456)$
to $(4)(56)$.

In the space
$\Xi\{5,6\}$, we have
$(2n+1,2n)\sim (1,0)$,
therefore, in $\Xi\{5,6\}$, we have  the constant sequence
$(1,0)$. Its limit is $(1,0)\in\Xi\{5,6\}$.
    \hfill $\square$

{\bf 6.11. Closures of strata.}
 The closure of a set $\Flat(\frA)$
is $\cup_\frB\Flat(\frB)$,
the union is given over all  refinements   $\frB$
of
the tree-partition $\frA$ (i.e., each element of $\frA$
is an element of $\frB$).

{\bf 6.12. Closure of the Weyl chamber
in the permutoassociahedron.}
Consider the Weyl
chamber $\Lambda_n^+=\Lambda^+(\bJ)$, i.e., the set of vectors
$\phi_1\ge\phi_2\ge\dots\ge\phi_n$
defined up to addition of a vector $(t,t,t,\dots)$.
Let us describe its closure $ \Lambda_n^\Pass$ in $\Pass_n$.

a) {\it Formal description.}
Denote $[j,k]=\{j,j+1,j+2,\dots,k\}\subset\bJ$.


We have the obvious projection $\Lambda_n^+\to\Lambda^+[j,k]$
and hence we have diagonal embedding
$$
\Lambda^+_n\,\,\mapsto \prod\limits_{1\le j< k\le n}
  \Lambda^+[j,k]\subset\prod\limits_{1\le j< k\le n}
\overline{\Lambda^+[j,k]}
$$
The set $\Lambda_n^\Pass$ is the closure of
$\Lambda_n^+$ in $\prod\limits_{1\le j< k\le n}
\overline{\Lambda^+[j,k]}$.

\smallskip

b) {\it List of strata.}
Strata are enumerated by segmental tree-partitions $\frA$
of $\bJ$.
A stratum  has the form
\begin{equation}
\prod\limits_{[k,l]\in\frA} U_\frA([k,l])
\end{equation}
and the factors  $U_\frA([k,l])$ are described in the following way:

 --- If $[k,l]$ is an irreducible element of $\frA$,
then $U_\frA([k,l]):=\Lambda^+[k,l]$.

--- Let $[k,l]$ be reducible.
Denote its minimal decomposition by $\frc$.
Then
$U_\frA([k,l]):=\Lambda^+\bigl([k,l]/\frc\bigr)^\gen$.


{\bf 6.13. Stasheff associahedron.}
Consider the tree-partition
$\frA_0:(1)(2)\dots(n)$
 of $\bJ$ and the corresponding
open face of $\Lambda^\Pass_n$.
The {\it associahedron} $\Ass_n$ is its closure
in $\Pass_n$.
Strata of $\Ass_n$
are enumerated
 by {\it perfect} segmental
tree-partitions of $\bJ$; they are described
in the previous subsection.

\smallskip

{\bf D. Description of the karpelevich-hedron}

\smallskip

{\bf 6.14. Stratification of the karpelevich-hedron.}
Strata $\Flat(\frA)$ of the karpelevich-hedron are enumerated
by leveled tree-partitions
\begin{equation}
\frA:\,\, \fra_0,\fra_1,\dots,\fra_\tau
\end{equation}
of the set $\bJ$.
Each stratum has the form
\begin{equation}
\Flat(\frA)=\Flat^\Karp(\frA)=\prod\nolimits_{j=0}^{\tau} W_j[\frA]
,\end{equation}
where the factors $W_j[\frA]$ are described in the following way.

a) Let $j<\tau$.
Consider the quotient set
$\bJ/\fra_{j+1}$.
The partition $\fra_{j}$ induces a partition of
$\bJ/\fra_{j+1}$. We denote this partition by
$\fra_{j}/\fra_{j+1}$.
We assume
$$
W_j[\frA]:=\partial\Xi[\bJ/\fra_{j+1};\fra_{j}/\fra_{j+1}]^\gen
$$
On  Fig. 4, Fig.6, the set   $\bJ/\fra_{j}$
is the set of
 of edges coming to the dotted line from below.
The
 set
$\bJ/\fra_{j+1}$  is the set is
 of edges coming to the dotted line above.
The quotient-partition $ \fra_{j}/\fra_{j+1}$
is the partition of the set $\bJ/\fra_{j+1}$
into churn-staffs.

\begin{figure}
\6
\end{figure}

b) For $j=\tau$, we assume
$
W_{\tau}[\frA]:=\Xi(\bJ;\fra_\tau)
$.

The {\it karpelevich-hedron} is a disjoint union
\begin{equation}
\Karp_n=\bigcup\nolimits_\frA \Flat(\frA)
\end{equation}
given over all the leveled tree-partitions
of $\bJ$.

{\bf 6.15. Identification of  definitions 6.4 and 6.14
of karpelevich-hedron.}
For each leveled tree-partition
$\frA$ and each partition $\frb$ of $\bJ$,
we must construct a map
\begin{equation}
\Flat(\frA)\to \overline{\Xi[\bJ;\frb]}
\end{equation}

Consider the maximal $j$ such that
$\frb$ is a refinement of $\fra_j$.
Such $j$ exists since $\fra_0$ is the trivial partition.
We consider the projection (see (6.19))
\begin{equation}
\Flat(\frA)\to W_{j}[\frA]
.\end{equation}

\smallskip

A) For $j<\tau$,
the map (6.21) is the composition of the maps
$$
\Flat(\frA)\to W_{j}[\frA]
=
\partial\Xi[\bJ/\fra_{j+1};\fra_j/\fra_{j+1}]^\gen
\to \partial \Xi[\bJ;\fra_j] \to
\partial\Xi[\bJ;\frb]
,$$
the first map is the projection to a factor in (6.19),
the second  map is (6.9), the third map is (6.9.a).


\smallskip

B) Let $j=\tau$.
Then $\frb$ is a refinement of $\fra_\tau$
and we have the canonical map
\begin{equation}
\Flat(\frA)\to
W_{\tau+1}[\frA]=\Xi[\bJ;\fra_\tau]\to \Xi[\bJ;\frb]
;\end{equation}
the second map is (6.5.a).


{\bf 6.16. Convergence in the karpelevich-hedron.}
Consider a sequence $x_j\in\Xi_n$.
Beginning of the definition of the convergence
is the same as in 6.10 until formula (6.15).

Then we obtain a sequence $x_j^s$ in each $\Xi(I_s)$,
or equivalently, a sequence in $\Xi[\bJ;\frp]$.
Our next condition is: the sequence $x_j$ converges in
$$\overline{\Xi[\bJ;\frp]}\simeq \overline{ \prod\Xi(I_s)}.$$
If the limit is contained in $\prod\Xi(I_s)$,
then it is the limit in $\Karp_n$. Otherwise,
let $t(\mu_1,\dots,\mu_n)$ be the limit ray.
 For $k$, $l\in\bJ$,
we say $k\sim l$ iff $k$, $l$ lie in one $I_s$, and
$\mu_k=\mu_l$.

Thus, we obtain a subpartition of each
element $I_s$,
hence we obtain a new partition $\frq$ of the whole set $\bJ$.

Our next condition of convergence is:
the sequence $x_j^s$ converges in the space
$\overline{\Xi[\bJ;\frq]}$, etc., etc., etc.

{\sc Example.}
For  sequence (6.16),
 the corresponding leveled tree-partition is
\begin{align*}
\fra_0=(123456); \quad \fra_1=(123)(456);
\quad \fra_2=(1)(2)(3)(4)(56)
\end{align*}

Indeed,
 we obtain the limit ray $t(1,1,1,0,0,0)\in\Xi_6$,
this gives the partition  $\fra_1=(123)(456)$.

In $\Xi[\bJ;\fra_1]\simeq \Xi\{1,2,3\}\times\Xi\{4,5,6\}$
we have the sequence
$$\bigl\{(n^3+2n,n^3+n,n^3) \times
  (3n,2n+1,2n)\bigr\}
\sim
\bigl\{(2n,n,0)\times (n,1,0)\bigr\}
 $$
Its limit ray  in  $\partial\Xi[\bJ;\fra_1]$
is  $t\bigl\{(2,1,0)\times (1,0,0)\bigr\}$.
This gives the partition    $ \fra_2=(1)(2)(3)(4)(56)$

Next,
$\Xi[\bJ;\fra_2]\simeq
\Xi\{1\}\times\Xi\{2\}\times\Xi\{3\}\times\Xi\{4\}
\times\Xi\{5,6\} \simeq \Xi\{5,6\}$,
In $\Xi\{5,6\}$, we have
$(2n+1, 2n)\simeq (1,0)$,
it is a constant sequence. Its limit is
 the point $(1,0)\in\Xi\{5,6\}$.
\hfill $\square$

{\bf 6.17.  Closures of strata.}
Let $\frA:\fra_0,\dots,\fra_p$ and $\frB:\frb_0,\dots,\frb_q$
be leveled tree-partitions. This say that $\frB$
is a refinement of  $\frA$,
if
each
partition $\fra_j$ is contained in the list
$\frb_0,\frb_2,\dots$.

The closure of the face
$\Flat(\frA)$ is $\cup\Flat(\frB)$ over all the
refinements $\frB$ of   $\frA$.

{\bf 6.18. Closure of Weyl chamber in the
karpelevich-hedron,}
see also \cite{GJT2}. Now let us describe the closure
$\Lambda^\Karp_n$ of the
positive Weyl chamber $\Lambda^+_n$
in the karpelevich-hedron.

{\it Abstract description.}
Let $\fra$ ranges in segmental partitions
of $\bJ$. Consider the natural map
$\Lambda^+_n\subset\Xi_n\to\Xi[\bJ;\fra]$ and the corresponding
diagonal map
$$
\Lambda^+_n\to\prod\nolimits_\fra\Xi[\bJ;\fra]
\subset
\prod\nolimits_\fra\overline{\Xi[\bJ;\fra]}
.$$
The set   $\Lambda^\Karp_n$
coincides with the closure of the image of
$\Lambda^+_n$ in
$\prod\nolimits_\fra\overline{\Xi[\bJ;\fra]}$

{\it Stratification.} Strata are enumerated
by segmental leveled tree-partitions
$\frA:\fra_0,\dots,\fra_p$ of $\bJ$.
Each stratum is the product
\begin{equation}
\prod_{j=0}^p Y_k[\frA]
,\end{equation}
where the factors have the following form

--- If $k<p$,  then
$Y_k[\frA]=\Lambda^+[\bJ/\fra_{j+1}; \fra_{j}/\fra_{j+1}]^\gen$.
In other words, we consider collections of
real numbers
$\theta(\mu)$, where $\mu$ ranges in   edges coming above
to a dotted line, see Fig.4, Fig.6. These numbers
are strictly increasing in each churn-staff
(if we move to right along the dotted
line), and they are defined
 up to a common positive factor and
an addition of a function that is constant on each
churn-staff.

--- $Y_p[\frA]=\prod_{[k,l]\in\fra_p} \Lambda^+[k,l]$.

{\bf 6.19. Map $\Karp_n\to\Pass_n$.}  Now we describe
the map $\pi:\Karp_n\to \Pass_n$
defined in 6.4.   Fix the notation of 6.8 and 6.14.
 Let $\frA$ be a leveled tree-partition
of $\bJ$. Let $\frA^\downarrow$ be the corresponding
tree-partition.

First, $\pi(\Flat^\Karp(\frA))=\Flat^\Pass(\frA^\downarrow)$.

Consider a partition $\fra_j$ lying in the leveled
tree-partition $\frA$.
Let $K_\alpha^{(j)}$ be its elements.
It is sufficient to describe the map
$$W_{j}[\frA]\to
\prod_\alpha U_{\frA^\downarrow}(K_\alpha^{(j)})
$$
(compare (6.12) and (6.19)). Consider two cases.

--- Let $j=\tau$. Then
$W_{\tau}\simeq \Xi[\bJ;\fra_\tau]$   coincides with
$\prod_{K\in\fra_\tau} \Xi(K)$,
and our map is the identical map.

--- For $j<\tau$, the map is described in Fig.6.

\begin{figure}
\7
\end{figure}

{\bf 6.20. Picture.}
{\it Karpelevich polyhedron.}
 $\Karp_4$ is a 3-dimensional polyhedron.
We can imagine surface of the polyhedron as a picture on
a sphere
(or
on the cube from Fig. 3.b.
We have 24 triangles on the surface of cube (sphere),
one of these triangles $PQT$ is drawn on Fig. 7.
We present division of this triangle into faces.

One of the 2-faces $ABCDEF$
of $\Karp_4$ is completely contained in the triangle $PQT$,
other faces have intersections with adjacent triangles.
The faces $\{\dots C'CDD'\dots\}$ and $\{\dots E'EFF'\dots\}$
are 12-gons.
The list of all faces (2-faces, edges, vertices)
 having intersection with $PQT$ is presented on Fig. 8.

{\it Permutoassociahedron.}
For obtaining the permutoassociahedron
$\Pass_4$ from the karpelevich-hedron $\Karp_4$,
it is sufficient to contract the edge $AB$ on Fig. 7
and 23 corresponding edges in other Weyl chambers.


\begin{figure}
\8
\end{figure}

        \smallskip

{\bf E. Root language}

\smallskip

{\bf 6.21. Permutoassociahedrons associated with  root systems.}
 Consider an {\it irreducible} root system $\Delta$
in a linear space $V(\Delta)$.
For each {\it irreducible} root subsystem
$\Gamma\subset\Delta$
consider its linear span $V(\Gamma)$,
and the corresponding compactification
$\overline{V(\Gamma)}$.
The sphere $\partial V(\Gamma)$
at infinity has the natural structure
of a {\it simplicial} complex.

Consider the orthogonal
projection
$\pi_\Gamma: V(\Delta)\to V(\Gamma)$
and the diagonal   embedding
$$
V(\Delta)\to
\prod_{\Gamma\subset\Delta} V(\Gamma)
\subset
\prod_{\Gamma\subset\Delta}
\overline {V(\Gamma)}
$$

The permutoassociahedron $\Pass(\Delta)$
is the closure of $V(\Delta)$ in
$\prod\overline {V(\Gamma)}$.

{\bf 6.22. Karpelevich-hedrons
 associated with  root systems.}
The definition of the karpelevich-hedron is the same,
we only omit two times the term {\it irreducible}
from the definition (and replace 'simplicial'
by 'polyhedral').

\smallskip

\medskip

{\bf 7. Existence of associahedral and Karpelevich boundaries}

\medskip

  \addtocounter{sec}{1}
 \setcounter{equation}{0}

\nopagebreak

{\bf 7.1. Hybrids of compactifications.}
Let $A$ be a metrizable space. Let $X$, $Y$ be compact
metrizable spaces
and
$\xi:A\to X$, $ \upsilon:A\to Y$
be continuous maps; assume
that the images of $A$ in $X$ and $Y$ are dense.

Then we have the diagonal map
$A\to X\times Y$
given by
$a\mapsto (\xi(a),\upsilon(a))$.
Consider the closure $Z\subset X\times Y$
of the image of $A$.
We say that $Z$ is a {\it hybrid} of compactifications
$X$ and $Y$.

{\bf 7.2. Velocity compactifications.}
Consider the positive Weyl chamber $\Lambda^+_n$ described
in 6.1. Let $\Lambda^\boxplus_n$ be a compact space
containing $\Lambda^+_n$ as a dense open subset.
Denote
$
\partial\Lambda^\boxplus_n=
\Lambda^\boxplus_n\setminus
\Lambda^+_n
$.

We define a structure of a compact space on the disjoint union
$\PE_n^\boxplus:=
\PE_n\cup \partial\Lambda^\boxplus_n
$.
Let $X^{(j)}$ be a sequence in $\PE_n$.
Let
$ \Phi^{(j)}: \phi_1^{(j)}\ge  \phi_1^{(j)}\ge\dots\ge  \phi_n^{(j)}$
be the eigenvalues of $X^{(j)}$.
We say that the sequence $X^{(j)}$ converges to
a point $\Psi\in \partial\Lambda^\boxplus_n$
if the sequence  $\Phi^{(j)}\in \Lambda^+_n$
converges to $\Psi$.

We say that $\PE_n^\boxplus$
 is a {\it velocity compactification} of $\PE_n$.
This construction is an analog of one-point compactification
of a locally compact space.

{\bf 7.3. Example. Martin boundary.}
The geometric object described below appears as the solution of
the problem of Martin boundary for symmetric spaces,  see
\cite{Olsh}, \cite{GJT2}.

Let $\Lambda^\boxplus_n$ be
$
\overline\Lambda^+_n=\Lambda^+_n\cup\partial\Lambda^+_n
$
defined in 6.2
The {\it Martin compactification} of $\PE_n$
is the hybrid of the velocity compactification
associated with
$\overline\Lambda^+_n$ and the Satake--Furstenberg compactification.

It is easy to describe it explicitly.
A point of the Martin compactification is a
following collection of data.

a) Subset
$I=\{0,i_1,\dots, i_{k},n\}
\subset \{0,1,\dots, n\}$.

b) A point of $\Delta_\circ(I)$, see 3.1.

c) A point of $\S(I)$, see 2.1.

{\bf 7.4. Construction of
the associahedral and Karpelevich boundaries.}
Consider  the completions
$\Lambda^\Pass_n$, $\Lambda^\Karp_n$ of the Weyl  chamber
described in 6.12 and 6.18.
Consider
the  associated velocity compactifications of $\PE_n$.
The associahedral and Karpelevich compactifications  of $\PE_n$
are the  hybrids of these velocity compactifications
with the Satake--Furstenberg compactification.

{\bf 7.5. Stratification of  associahedral compactification.}
A point of the  associahedral compactification
of $\PE_n$ is the following
collection
of data A-C

A. A segmental tree-partition $\frA$  of the set $\bJ$.
Denote by
$[1,i_1]$,
$[i_1+1,i_2]$,\dots, $[i_s+1,n]$
its irreducible elements.

B. A point of the
stratum (6.17) of $\Lambda_n^\Pass$.
We emphasis, that
the product (6.17) contains the factors
$\prod_k \Lambda^+([i_{k-1}+1,i_{k+1}])$.
A point of a factor
$\Lambda^+([i_{k-1}+1,i_{k+1}])$
is a collection  of numbers
\begin{equation}
 \psi_{{i_k}+1}\ge\psi_{{i_k}+2}\ge\dots \ge \psi_{i_{k+1}}
\end{equation}
defined up to addition of a common  constant.

C. A point of the set $\S\{i_1,\dots,i_s\}$, i.e. a flag
$\R^n=W_0\supset W_1\supset\dots\supset W_{i_s+1}=0$
and an ellipsoid $Q_k$ in each subquotient
$W_{k-1}/W_k$ of the flag.
These ellipsoids are not arbitrary.
Our additional requirement is: for each $k$
the  principal semiaxes of $Q_k$
are $e^{\psi_{i_{k-1}}}$, \dots, $e^{\psi_{i_{k}}}$, where
$\psi_{...}$ are already defined (7.1).

{\bf 7.6. Stratification of the Karpelevich compactification.}
A point of the Karpelevich compactification
of $\PE_n$ is the following
collection
of data A-C

A. A segmental leveled tree-partition $\frA$  of the set $\bJ$.
Denote by
$[1,i_1]$, $ [i_1+1,i_2]$, \dots, $[i_s+1,n]$
its irreducible elements.

B. A point of the
stratum (6.24) of $\Lambda_n^\Karp$.
We emphasis, that
the product (6.24) contains the factor
$Y_p(\frA)
=\prod_k \Lambda^+([i_{k-1}+1,i_{k+1}])$.
A point of $k$-th factor of the last product
is a collection (7.1).

C. A point of the set $\S\{i_1,\dots,i_s\}$.
The ellipsoids $Q_k$ in the subquotients
$W_{k-1}/W_k$ of the flag
 are not arbitrary.
Our additional requirement is: for each $k$
the  principal semiaxes of $Q_k$
are $e^{\psi_{i_{k-1}}}$, \dots, $e^{\psi_{i_{k}}}$, where
$\psi_{...}$ are already defined (7.1).

\medskip

{\bf 8. Sea urchin}

\medskip

  \addtocounter{sec}{1}
 \setcounter{equation}{0}

\nopagebreak

A simple calculation shows that the dimension of
the space of all null-pencils is larger than
$\dim \PE_n$.
Nevertheless it is possible to define
a boundary related to null-pencils.

{\bf 8.1. Definition.} We define the sea urchin boundary
as the set of limits of null-pencils
whose
 velocities
$
(\phi_1,\dots,\phi_n)
$
consist of integer numbers.

Explicit description of the sea urchin
can be easily obtained from  3.9.

{\bf 8.2. Limits of meromorphic curves in sea urchin.}
Let $z$ ranges in a small interval $(0,\epsilon)$.
We say, that a map $z\mapsto X(z)$ is a meromorphic curve
in $\E_n$ if each matrix $x_{ij}$ element admits
a Laurent decomposition
$$
x_{ij}(z)= z^{-k} (a_0+a_1 z +a_2 z^2+\dots)
.$$

{\sc Lemma 8.1.}
{\it Each meromorphic curve $X(z)\in\E_n$
admits a  representation
\begin{equation}
X(z)=g(z)
\begin{pmatrix}
z^{-k_1}&0&\dots\\
0&z^{-k_2}&\dots\\
\vdots&\vdots&\ddots
\end{pmatrix}
g(z)^\top
\end{equation}
where  $k_1\ge k_2\ge\dots$,
the function $g(z)$ is analytic in a neighborhood of 0,
and $g(0)$ is invertible.}

(We reduce the positive definite quadratic form
$X(z)$ to a sum of squares in the usual way).

For the curve (8.1),
consider the geodesic
\begin{equation}
\gamma(t)   =
g(0)
\begin{pmatrix}
e^{k_1 t}&0&\dots\\
0&e^{k_2 t}&\dots\\
\vdots&\vdots&\ddots
\end{pmatrix}
g(0)^\top
,\end{equation}
we substitute $z=e^{-t}$ to the middle factor in (8.1),
and $z=0$ to the first and last factors.
We define the limit of
the meromorphic curve (8.1) in the sea urchin
as the limit of the geodesic (8.2).

{\sc Lemma 8.2.}
a) {\it A limit of a meromorphic curve does not depend on
choice of the representation} (8.2).

b) {\it The limit of a meromorphic curve does not
depend
on a parametrization of the curve.}

{\sc Remark.} The sea urchin is not a compact
space in the usual sense.
For instance, the sequence
$
X_k=\begin{pmatrix}
   e^k&0\\
   0&k
\end{pmatrix}
$
has no limit (and no limit points)  in the sea urchin.

{\bf 8.3. Projective compactifications.}
Consider a polynomial representation
$\rho_{\bold m}$ of $\GL(n,\R)$ with
a highest weight
${\bold m}:\,\,m\ge m_2\ge\dots\ge m_n$, where $m\in {\Bbb Z}
.$
It is well known, that
 the representation $\rho_{\bold m}$ contains a nonzero
$\O(n)$-invariant vector iff
all the numbers $m_j$ are even.
In this case, an $\O(n)$-invariant vector $q_{\bold m}$
 is unique
up to a scalar factor.

\smallskip

Consider a direct sum $\theta$ of several
 representations
$\rho_{{\bold m}^\tau}$ of $\GL(n,\R)$ with even
signatures ${\bold m}^\tau$. Denote by $H$ the space
of the representation
$\theta$. Denote by $h$ the $\O(n)$-invariant vector
$\oplus q_{{\bold m}^\tau}$
in $H$.

Consider the projective space $\P H$ and the
$\GL(n,\R)$-orbit $\cal O$ of the vector $h$
in $\P H$.
The projective compactification $[\PE_n]_\theta$
of $\PE_n$ is the closure of the orbit
$\cal O$ in $\P H$.

{\bf 8.5. Universality of the sea urchin.}
Obviously, each meromorphic curve $X(z)$ has a limit
in each projective compactification.

{\sc Theorem 8.3.}
a) {\it For each point $A$ of each projective compactification
$[\E_n]_\theta$,
there exists
a meromorphic curve $X(z)\in\E_n$, whose
limit in  $[\E_n]_\theta$ is $A$.}

b) {\it If the limits of two meromorphic curves
$X_1(z)$ and $X_2(z)$ in the sea urchin coincide, then
their limits in any projective compactification coincide.}

c) {\it If the limits of $X_1(z)$ and $X_2(z)$
in the sea urchin are different, then their limits in
some projective compactification $[\E_n]_\theta$
are different.}

Thus, for each projective compactification
$[\E_n]_\theta$, we have the canonical surjective
map from the sea urchin to   $[\E_n]_\theta$;
an explicit variant of this construction is contained in
\cite{Ner7}.

 \smallskip

{\bf Acknowledgments.}
I am grateful to C.De Concini for discussions
of complete symmetric varieties and of
the universalization problem. I thank M.Kapranov,
C.Kapoudjian and M.A.Olshanetsky for interesting
discussions. I also thank J.Stasheff
who identified a polyhedron from \cite{Ner4}
with permutoassociahedron and to L.Ji
who mentioned me a mistake in \cite{Ner4}
(difference between the associahedral
and Karpelevich  boundaries), and W.Fulton
for the reference \cite{Ulya}. These
notes are partially based on my lectures
given in Winter school in Crni (Chehia)
in January 2002. I thank its organizers.
 The text of the work was prepared during
my visit to the University of Michigan.
I thank the administrators and
  W.Fulton for hospitality

\sf\small Math. Phys. Group,  Institute of Theoretical
and   Experimental Physics,

 Bolshaya Cheremushkinskaya, 25, Moscow 117 259, Russia

\tt neretin@mccme.ru, neretin@gate.itep.ru

\end{document}